\documentclass[12pt]{article}
\usepackage[letterpaper, left=2.5cm, right=2.5cm, top=1.5cm, bottom=1.5cm, headheight=7mm, headsep=0.8cm, footskip=1.5cm, includeheadfoot]{geometry}
\usepackage[utf8]{inputenc}
\usepackage{amsmath}
\usepackage{authblk}
\usepackage{amssymb}
\usepackage{bm}
\usepackage{amsthm}
\usepackage{mathtools}
\usepackage{algorithm}
\usepackage{algpseudocode}
\usepackage{color}
\usepackage{graphicx}
\usepackage{subcaption}
\usepackage{MnSymbol}
\usepackage{hyperref}
\hypersetup{
    unicode=false,  
    pdftoolbar=true,
    pdfmenubar=true,
    pdffitwindow=false,
    pdfstartview={FitH},
    pdfnewwindow=true,
    colorlinks=true,
    linkcolor=blue,
    citecolor=blue,
    filecolor=blue,
    urlcolor=blue
}

\usepackage{tikz}
\usepackage{tikz-3dplot} 
\usetikzlibrary{intersections,patterns.meta}
\usetikzlibrary{positioning}
\usetikzlibrary{calc}
\usetikzlibrary{decorations.pathmorphing}

\newtheorem{theorem}{Theorem}
\newtheorem{proposition}{Proposition}
\newtheorem{lemma}{Lemma}
\newtheorem{definition}{Definition}
\newtheorem{remark}{Remark}
\newtheorem{corollary}{Corollary}

\theoremstyle{definition}

\DeclareMathOperator*{\ext}{ext}

\title{Nearly-periodic maps and geometric integration of noncanonical Hamiltonian systems}
\author[1]{J. W. Burby}
\author[2]{E. Hirvijoki}
\author[3]{M. Leok}
\affil[1]{Los Alamos National Laboratory, Los Alamos, NM 87545, USA}
\affil[2]{Department of Mechanical Engineering, Aalto University, Finland}
\affil[3]{Department of Mathematics, University of California, San Diego, 9500 Gilman Drive, La Jolla, CA 92093-0112, USA}
\date{March 2021}

\begin{document}

\maketitle

\begin{abstract}
M.~Kruskal showed that each continuous-time nearly-periodic dynamical system admits a formal $U(1)$ symmetry, generated by the so-called roto-rate. When the nearly-periodic system is also Hamiltonian, Noether's theorem implies the existence of a corresponding adiabatic invariant. We develop a discrete-time analoue of Kruskal's theory. Nearly-periodic maps are defined as parameter-dependent diffeomorphisms that limit to rotations along a $U(1)$-action. When the limiting rotation is non-resonant, these maps admit formal $U(1)$ symmetries to all orders in perturbation theory. For Hamiltonian nearly-periodic maps on exact presymplectic manifolds, we prove that the formal $U(1)$ symmetry gives rise to a discrete-time adiabatic invariant using a discrete-time extension of Noether's theorem. When the unperturbed $U(1)$-orbits are contractible, we also find a discrete-time adiabatic invariant for mappings that are merely presymplectic, rather than Hamiltonian. As an application of the theory, we use it to develop a novel technique for geometric integration of non-canonical Hamiltonian systems on exact symplectic manifolds.
\end{abstract}

\tableofcontents

\section{Introduction}
%journals: 1. {\color{blue}Journal of Nonlinear Science}, 2. Nonlinearity
%, SIAM applied dynamical systems, for theory
%for applications: JCP

A continuous-time dynamical system with vector parameter $\gamma$ is {nearly periodic} if all of its trajectories are periodic with nowhere-vanishing angular frequency in the limit $\gamma\rightarrow 0$. Examples from physics include charged particle dynamics in a strong magnetic field, the weakly-relativistic Dirac equation, and any mechanical system subject to a high-frequency, time-periodic force. In the broader context of multi-scale dynamical systems, nearly-periodic systems play a special role because they display perhaps the simplest possible non-dissipative short-timescale dynamics. They therefore provide a useful proving ground for analytical and numerical methods aimed at more complex multi-scale models.

In a seminal paper \cite{Kruskal_1962}, Kruskal deduced the basic asymptotic properties of continuous-time nearly-periodic systems. In general, each such system admits a formal $U(1)$ Lie symmetry whose infinitesimal generator $R_\gamma$ is known as the {roto-rate}. In the Hamiltonian setting, existence of the roto-rate implies existence of an all-orders adiabatic invariant $\mu_\gamma$ by way of Noether’s theorem. General expressions for $\mu_\gamma$ may be found in \cite{Burby_Squire_2020}. Recently \cite{Burby_Hirvijoki_2021_JMP}, we extended Kruskal’s analysis by proving that the (formal) set of fixed points for the roto-rate is an elliptic almost invariant slow manifold. Moreover, in the Hamiltonian case, we demonstrated that normal stability of the slow manifold is mediated by Kruskal’s adiabatic invariant.

The purpose of this article is to introduce discrete-time analogues of continuous-time nearly-periodic systems that we call {nearly-periodic maps}. These objects can be motivated as follows. A nearly-periodic system characteristically displays limiting short-timescale dynamics that ergodically cover circles in phase space. This ergodicity is ultimately what gives rise to Kruskal's roto-rate and, in the presence of Hamiltonian structure, adiabatic invariance. It is therefore sensible to regard parameter-dependent maps whose limiting iterations ergodically cover circles as discrete-time analogues of nearly-periodic systems. Ergodicity requires that the rotation angle associated with each circle be an irrational multiple of $2\pi$. In principle these rotation angles could vary from circle to circle, but smoothness removes this freedom, and imposes a common rotation angle across circles. Nearly-periodic maps are defined by limiting iterations that rotate a family of circles foliating phase space by a common rotation angle. Such a map is resonant or non-resonant when the rotation angle is a rational or irrational multiple of $2\pi$, respectively. The preceding remarks suggest that non-resonant nearly-perioidic maps should share important features with continuous-time nearly-periodic systems.

We will show that non-resonant nearly-periodic maps always admit formal $U(1)$ symmetries by modifying Kruskal's construction of a normal form for the roto-rate. Thus, non-resonant nearly-periodic maps formally reduce to mappings on the space of $U(1)$-orbits, corresponding to elimination of a single dimension in phase space. In the Hamiltonian setting, we will establish a discrete-time analogue of Noether’s theorem that will allow us to construct all-orders adiabatic invariants for non-resonant nearly-periodic maps. In contrast to the continuous-time case, there may be topological obstructions to the Noether theorem based construction. Nevertheless, assuming (a) existence of a fixed point for formal $U(1)$ symmetry, or (b) existence of a time-dependent Hamiltonian suspension for the nearly-periodic map, these topological obstructions disappear. When an adiabatic invariant does exist, the phase-space dimension is formally reduced by two instead of one.

We anticipate that non-resonant nearly-periodic maps will have important applications to numerical integration of nearly-periodic systems. While development of integrators for such systems is straightforward when the numerical timestep $h$ resolves the short-timescale dynamics, considerably more care is required when ``stepping over" the period of limiting oscillations. 
%A time-tested approach to addressing this problem is to first analytically derive evolution equations on the space of $U(1)$-orbits, which are automatically free of the short oscillation timescale, and the apply standard integration methods. Downsides of this standard approach include (a) the reduced space may be geometrically more complex than the unreduced space, as happens for charged particles in a strong magnetic field, and (b) trajectories that move into regions where the limiting angular frequency is small cannot be handled at all. 
One approach would be to design an integrator on the unreduced space that is constrained to be a non-resonant nearly-periodic map. Although such an integrator would not accurately resolve the phase of short-scale oscillations when taking large timesteps, it would automatically possess an all-orders reduction to the space of $U(1)$ orbits. By designing the reduced map to discretize the continuous-time reduced dynamics, the slow component of the continuous-time dynamics could be accurately resolved without directly simulating the reduced dynamical variables. This opens the door to a type of asymptotic-preserving integrator capable of seamlessly transitioning between large- and small-timestep regimes, generalizing those proposed in \cite{Ricketson_2020,Xiao_Qin_2021} for magnetized charged particle dynamics. Moreover, in the Hamiltonian case, the integrator would automatically enjoy an all-orders adiabatic invariant close to the continuous-time invariant. Such a capability would complement previous results on short-timestep adiabatic invariants for variational integrators \cite{Hairer_Lubich_2020}. We provide a proof-of-principle demonstration of these ideas in Section \ref{hidden_gravity_section}

Aside from serving as integrators for nearly-periodic systems, nearly-periodic maps may also be used as tools for structure-preserving simulation of general Hamiltonian systems on exact symplectic manifolds. (See \cite{Abraham_2008,Marsden_Ratiu_1999} for the foundations of Hamiltonian mechanics on symplectic manifolds.) The basic idea is to first embed the original Hamiltonian system as an approximate invariant manifold inside of a larger nearly-periodic Hamiltonian system, as discussed in \cite{Burby_Hirvijoki_2021_JMP}. Then it is possible to construct a symplectic nearly-periodic map that integrates the larger system while preserving the approximate invariant manifold. Discrete-time adiabatic invariance ensures that the approximate invariant manifold enjoys long-term normal stability, which is tantamount to the integrator providing a persistent approximation of the original system's dynamics. We describe and analyze this construction in Section \ref{sec:symplectic_lorentz}. In Section \ref{reduced_gc}, we apply the general theory to the non-canonical Hamiltonian dynamics of a charged particle's guiding center \cite{Northrop_1963,Littlejohn_1981,Littlejohn_1983} in a magnetic field of the form $\bm{B} = B(x,y)\,\bm{e}_z$ \cite{Littlejohn_1979}.

The remainder of this article is organized as follows. We review Kruskal's theory of nearly-periodic systems using modern terminology in Section \ref{nps_background}. Then we develop the general theory of nearly-periodic maps in Section \ref{npm_general_sec}, including their special properties in the symplectic case, and their ability to serve as geometric integrators for Hamiltonian systems on exact symplectic manifolds. Wherever possible, proofs of general properties of nearly-periodic maps parallel Kruskal's arguments from the continuous-time setting. Section \ref{examples_section} contains a pair of interesting applications of nearly-periodic map technology. Finally, Section \ref{discussion_section} provides additional review and context for this work.

\subsection{Notational conventions}
%{\color{red}fold this into introduction}
In this article, smooth shall always mean $C^\infty$, and $\Gamma$ will always denote a vector space. We reserve the symbol $M$ for a smooth manifold equipped with a smooth auxilliary Riemannian metric $g$. We say $f_\gamma:M_1\rightarrow M_2$, $\gamma\in\Gamma$, is a smooth $\gamma$-dependent mapping between manifolds $M_1,M_2$ when the mapping $M_1\times\mathbb{R}\rightarrow M_2:(m,\gamma)\mapsto f_\gamma(m)$ is smooth. Similarly, $\bm{T}_\gamma$ is a smooth $\gamma$-dependent tensor field on $M$ when (a) $\bm{T}_\gamma(m)$ is an element of the tensor algebra $\mathcal{T}_m(M)$ at $m$ for each $m\in M$  and $\gamma\in\Gamma$, and (b) $\bm{T}_\gamma$ is a smooth $\gamma$-dependent mapping between the manifolds $M$ and $\mathcal{T}(M)=\cup_{m\in M}\mathcal{T}_m(M)$. 

The symbol $X_\gamma$ will always denote a smooth $\gamma$-dependent vector field on $M$. If $\bm{T}_\gamma$ is a smooth $\gamma$-dependent section of either $TM\otimes TM$ or $T^*M\otimes T^*M$ then $\widehat{\bm{T}}_\gamma$ is the corresponding smooth $\gamma$-dependent bundle map $T^*M\rightarrow TM:\alpha\mapsto \iota_\alpha\bm{T}_\gamma$, or $TM\rightarrow T^*M:X\mapsto \iota_X\bm{T}_\gamma$, respectively. Note that if $\Omega$ is a symplectic form on $M$ with associated Poisson bivector $\mathcal{J}$ then $\widehat{\Omega}^{-1} = -\widehat{\mathcal{J}}$.
% \begin{align}
%     \iota_{X}\Omega &= \mathbf{d}H\\
%     \widehat{\Omega}\,X & = \mathbf{d}H\\
%     X & = \widehat{\Omega}^{-1}\,\mathbf{d}H\\
%     X & = \widehat{\mathcal{J}}\,\mathbf{d}H\\
%     \mathcal{J}(df,dg) &= \Omega(X_f,X_g) = \mathbf{d}f(\widehat{\Omega}^{-1}\,\mathbf{d}g)\\
%     -\iota_{dg}\mathcal{J} = \widehat{\Omega}^{-1}\,dg
% \end{align}

\section{Kruskal's theory of nearly-periodic systems\label{nps_background}}

In 1962, Kruskal presented an asymptotic theory \cite{Kruskal_1962} of averaging for dynamical systems whose trajectories are all periodic to leading order. Nowadays, Kruskal's method is termed one-phase averaging \cite{Lochak_1993}, which suggests a contrast with the multi-phase averaging methods underlying, e.g., Kolmogorov--Arnol'd--Moser (KAM) theory. Since this theory provides a model for the results in this article, we review its main ingredients here. In this section only, and merely for simplicity's sake, we make the restriction $\Gamma=\mathbb{R}$.

\begin{definition}
A \textbf{nearly-periodic system} on a manifold $M$ is a smooth $\gamma$-dependent vector field $X_\gamma$ on $M$ such that $X_0 = \omega_0\,R_0$, where
\begin{itemize}
    \item $\omega_0:M\rightarrow\mathbb{R}$ is strictly positive
    \item $R_0$ is the infinitesimal generator for a circle action $\Phi_\theta:M\rightarrow M$, $\theta\in U(1)$.
    \item $\mathcal{L}_{R_0}\omega_0 = 0$.
\end{itemize}
The vector field $R_0$ is called the \textbf{limiting roto-rate}, and $\omega_0$ is the \textbf{limiting angular frequency}.%, and the set $S_0 = \{s\in M\mid R_0(s) = 0\}$ is called the \textbf{limiting slow manifold}.
\end{definition}
\begin{remark}
In addition to requiring that $\omega_0$ is sign-definite, Kruskal assumed that $R_0$ is nowhere vanishing. However, this assumption is not essential for one-phase averaging to work. It is enough to require that $\omega_0$ vanishes nowhere. This is an important restriction to lift since many interesting circle actions have fixed points. \end{remark}

Kruskal's theory applies to both Hamiltonian and non-Hamiltonian systems. In the Hamiltonian setting, it leads to stronger conclusions. A general class of Hamiltonian systems for which the theory works nicely may be defined as follows.

\begin{definition}\label{nearly_periodic_hamiltonian_system_def}
Let $(M,\Omega_\gamma)$ be a manifold equipped with a smooth $\gamma$-dependent presymplectic form $\Omega_\gamma$. Assume there is a smooth $\gamma$-dependent $1$-form $\vartheta_\gamma$ such that $\Omega_\gamma = - \mathbf{d}\vartheta_\gamma$. A \textbf{nearly-periodic Hamiltonian system} on $(M,\Omega_\gamma)$ is a nearly-periodic system $X_\gamma$ on $M$ such that $\iota_{{X}_\gamma}\Omega_\gamma = \mathbf{d}H_\gamma$, for some smooth $\gamma$-dependent function $H_\gamma:M\rightarrow\mathbb{R}$. 
\end{definition}

Kruskal showed that all nearly-periodic systems admit an approximate $U(1)$-symmetry that is determined to leading order by the unperturbed periodic dynamics. He named the generator of this approximate symmetry the \emph{roto-rate}. In the Hamiltonian setting, he showed that both the dynamics and the Hamiltonian structure are $U(1)$-invariant to all orders in $\gamma$.

\begin{definition}
A \textbf{roto-rate} for a nearly-periodic system $X_\gamma$ on a manifold $M$ is a formal power series $R_\gamma = R_0 + \gamma\,R_1 + \gamma^2\,R_2 + \dots$ with vector field coefficients such that
\begin{itemize}
    \item $R_0$ is equal to the limiting roto-rate
    \item $\exp(2\pi \mathcal{L}_{R_\gamma}) = 1$
    \item $[X_\gamma,R_\gamma] = 0$,
\end{itemize}
where the second and third conditions are understood in the sense of formal power series.
\end{definition}

\begin{proposition}[Kruskal \cite{Kruskal_1962}]\label{existence_of_roto_rate}
Every nearly-periodic system admits a unique roto-rate $R_\gamma$. The roto-rate for a nearly-periodic Hamiltonian system on an exact presymplectic manifold $(M,\Omega_\gamma)$ satisfies $\mathcal{L}_{R_\gamma}\Omega_\gamma = 0$ in the sense of formal power series. 
\end{proposition}

\begin{corollary}\label{invariance_of_Hamiltonian}
The roto-rate $R_\gamma$ for a nearly-periodic Hamiltonian system $X_\gamma$ on an exact presymplectic manifold $(M,\Omega_\gamma)$ with Hamiltonian $H_\gamma$ satisfies $\mathcal{L}_{R_\gamma}H_\gamma = 0$.
\end{corollary}
\begin{proof}
Since $[R_\gamma,X_\gamma] = \mathcal{L}_{R_\gamma}X_\gamma = 0$ and $\mathcal{L}_{R_\gamma}\Omega_\gamma = 0$, we may apply the Lie derivative $\mathcal{L}_{R_\gamma}$ to Hamilton's equation $\iota_{X_\gamma}\Omega_\gamma = \mathbf{d}H_\gamma$ to obtain
\begin{align*}
\mathcal{L}_{R_\gamma}(\mathbf{d}H_\gamma) = \mathcal{L}_{R_\gamma}(\iota_{X_\gamma}\Omega_\gamma) = \iota_{\mathcal{L}_{R_\gamma} X_\gamma}\Omega_\gamma + \iota_{X_\gamma}(\mathcal{L}_{R_\gamma}\Omega_\gamma) = 0.
\end{align*}
Thus, $\mathcal{L}_{R_\gamma}H_\gamma$ is a constant function. By averaging over the $U(1)$-action we conclude that the constant must be zero.
\end{proof}

To prove Proposition \ref{existence_of_roto_rate}, Kruskal used a pair of technical results, each of which is interesting in its own right. The first establishes the existence of a non-unique normalizing transformation that asymptotically deforms the $U(1)$ action generated by $R_\gamma$ into the simpler $U(1)$-action generated by $R_0$. The second is a subtle bootstrapping argument that upgrades leading-order $U(1)$-invariance to all-orders $U(1)$-invariance for integral invariants. We state these results here for future reference.

%\todo{Why is $G_\gamma$ a $O(\gamma)$ formal power series as opposed to just a formal power series in $\gamma$? In particular, in the definition, it is not approximating anything yet.}

\begin{definition}
Let $G_\gamma =\gamma \,G_1+ \gamma^2\,G_2 + \dots $ be an $O(\gamma)$ (no constant term) formal power series whose coefficients are vector fields on a manifold $M$. 
The \textbf{Lie transform} with \textbf{generator} $G_\gamma$ is the formal power series $\exp(\mathcal{L}_{G_\gamma})$ whose coefficients are differential operators on the tensor algebra over $M$. 
\end{definition}

\begin{definition}\label{normalizing_transformation_sys_def}
A \textbf{normalizing transformation} for a nearly-periodic system $X_\gamma$ with roto-rate $R_\gamma$ is a Lie transform $\exp(\mathcal{L}_{G_\gamma})$ with generator $G_\gamma$ such that $R_\gamma = \exp(\mathcal{L}_{G_\gamma})R_0$.
\end{definition}

\begin{proposition}[Kruskal]\label{existence_of_normalizing_transformation}
Each nearly-periodic system admits a normalizing transformation.
\end{proposition}

\begin{proposition}\label{bootstrap_prop}
Let $\alpha_\gamma$ be a smooth $\gamma$-dependent differential form on a manifold $M$. Suppose $\alpha_\gamma$ is an absolute integral invariant for a $C^\infty$ nearly-periodic system $X_\gamma$ on $M$. If $\mathcal{L}_{R_0}\alpha_0 = 0$ then $\mathcal{L}_{R_\gamma}\alpha_\gamma = 0$, where $R_\gamma$ is the roto-rate for $X_\gamma$.
\end{proposition}
\begin{proof}
Integral invariance means $\mathcal{L}_{X_\gamma}\alpha_\gamma = 0$ for each $\gamma\in\Gamma$. By Applying $\mathcal{L}_{R_\gamma}$ to this relationship, and using $[R_\gamma,X_\gamma]=0$, we obtain $\mathcal{L}_{X_\gamma}\mathcal{L}_{R_\gamma}\alpha_\gamma = 0$. Now let $G_\gamma$ be the generator of a normalizing transformation for $X_\gamma$, and set $\overline{X}_\gamma = \exp(-\mathcal{L}_{G_\gamma})X_\gamma$, $\overline{\alpha}_\gamma = \exp(-\mathcal{L}_{G_\gamma})\alpha_\gamma$. We have $\mathcal{L}_{\overline{X}_\gamma}\mathcal{L}_{R_0}\overline{\alpha}_\gamma=0$. Since $\mathcal{L}_{R_0}\overline{\alpha}_\gamma = O(\gamma)$, the first-order consequence of the previous formula is $\mathcal{L}_{\omega_0\,R_0}\mathcal{L}_{R_0}\overline{\alpha}_1=0$, which can only be satisfied if $\mathcal{L}_{R_0}\overline{\alpha}_1$ is $R_0$-invariant. But since the $U(1)$-average of $\mathcal{L}_{R_0}\overline{\alpha}_1$ vanishes, we conclude  $\mathcal{L}_{R_0}\overline{\alpha}_1=0$. Repeating this argument gives $\mathcal{L}_{R_0}\overline{\alpha}_k=0$ for $k > 1$ as well. In other words $\mathcal{L}_{R_0}\overline{\alpha}_\gamma=0$ to all orders in $\gamma$, which is equivalent to the Theorem's claim.
\end{proof}

According to Noether's celebrated theorem, a Hamiltonian system that admits a continuous family of symmetries also admits a corresponding conserved quantity. Therefore one might expect that a Hamiltonian system with approximate symmetry must also have an approximate conservation law. Kruskal showed that this is indeed the case for nearly-periodic Hamiltonian systems, as the following generalization of his argument shows.

\begin{proposition}\label{existence_of_mu}
Let $X_\gamma$ be a nearly-periodic Hamiltonian system on the exact presymplectic manifold $(M,\Omega_\gamma)$. Let $R_\gamma$ be the associated roto-rate. There is a formal power series $\theta_\gamma = \theta_0 + \gamma\,\theta_1 + \dots$ with coefficients in $\Omega^1(M)$ such that $\Omega_\gamma = -\mathbf{d}\theta_\gamma$ and $\mathcal{L}_{R_\gamma}\theta_\gamma = 0$. Moreover, the formal power series $\mu_\gamma = \iota_{R_\gamma}\theta_\gamma$ is a constant of motion for $X_\gamma$ to all orders in perturbation theory. In other words,
\begin{align*}
\mathcal{L}_{X_\gamma}\mu_\gamma = 0,
\end{align*}
in the sense of formal power series.
\end{proposition}
\begin{proof}
To construct the $U(1)$-invariant primitive $\theta_\gamma$ we select an arbitrary primitive $\vartheta_\gamma$ for $\Omega_\gamma$ and set
\begin{align*}
\theta_\gamma = \frac{1}{2\pi}\int_0^{2\pi}\exp(\theta\mathcal{L}_{R_\gamma})\vartheta_\gamma\,d\theta.
\end{align*}
This formal power series satisfies $\mathcal{L}_{R_\gamma}\theta_\gamma=0$ because
\begin{align*}
\mathcal{L}_{R_\gamma}\theta_\gamma = \frac{1}{2\pi}\int_0^{2\pi}\frac{d}{d\theta}\exp(\theta\mathcal{L}_{R_\gamma})\vartheta_\gamma\,d\theta = 0.
\end{align*}
Moreover, since $\mathcal{L}_{R_\gamma}\Omega_\gamma = 0$ by Kruskal's Proposition \ref{existence_of_roto_rate}, we have
\begin{align*}
-\mathbf{d}\theta_\gamma = \frac{1}{2\pi}\int_0^{2\pi}\exp(\theta\mathcal{L}_{R_\gamma})\Omega_\gamma\,d\theta = \frac{1}{2\pi}\int_0^{2\pi}\Omega_\gamma\,d\theta = \Omega_\gamma,
\end{align*}
whence $\theta_\gamma$ is a primitive for $\Omega_\gamma$.

To establish all-orders time-independence of $\mu_\gamma = \iota_{R_\gamma}\theta_\gamma$, we apply Cartan's formula and Corollary \ref{invariance_of_Hamiltonian} according to
\begin{align*}
\mathcal{L}_{X_\gamma}\mu_\gamma = \iota_{X_\gamma}\mathbf{d}\iota_{R_\gamma}\theta_\gamma  =  - \iota_{R_\gamma}\iota_{X_\gamma}\Omega_\gamma = - \mathcal{L}_{R_\gamma}H_\gamma =  0. 
\end{align*}
\end{proof}

\begin{definition}\label{mu_defined}
The formal constant of motion $\mu_\gamma$ provided by Proposition \ref{existence_of_mu} is the \textbf{adiabatic invariant} associated with a nearly-periodic Hamiltonian system.
\end{definition}

\section{Nearly-periodic maps\label{npm_general_sec}}

Nearly-periodic maps are natural discrete-time analogues of nearly-periodic systems. The following provides a precise definition.
\begin{definition}[nearly-periodic map]
Let $\Gamma$ be a vector space. A \textbf{nearly-periodic map} on a manifold $M$ with parameter space $\Gamma$ is a smooth mapping $F:M\times \Gamma\rightarrow M$ such that $F_\gamma:M\rightarrow M:m\mapsto F(m,\gamma)$ has the following properties:
\begin{itemize}
\item $F_\gamma$ is a diffeomorphism for each $\gamma\in \Gamma$.
\item There exists a $U(1)$-action $\Phi_\theta:M\rightarrow M$ and a constant $\theta_0\in U(1)$ such that $F_0 = \Phi_{\theta_0}$.
\end{itemize}
We say $F$ is \textbf{resonant} if $\theta_0$ is a rational multiple of $2\pi$, otherwise $F$ is \textbf{non-resonant}. The infinitesimal generator of $\Phi_\theta$, $R_0$, is the \textbf{limiting roto-rate}.

\end{definition}

Let $X$ be a vector field on a manifold $M$ with time-$t$ flow map $\mathcal{F}_t$. A $U(1)$-action $\Phi_\theta$ is a symmetry for $X$ if $\mathcal{F}_t\circ \Phi_\theta = \Phi_\theta\circ \mathcal{F}_t$, for each $t\in\mathbb{R}$ and $\theta\in U(1)$. Differentiating this condition with respect to $\theta$ at the identity implies, and is implied by, $\mathcal{F}_t^*R = R$, where $R$ denotes the infinitesimal generator for the $U(1)$-action. Since we would like to think of nearly-periodic maps as playing the part of a nearly-periodic system's flow map, the latter characterization of symmetry allows us to naturally extend Kruskal's notion of roto-rate to our discrete-time setting.
\begin{definition}
A \textbf{roto-rate} for a nearly-periodic map $F$ is a formal power series $R_\gamma = R_0 + R_1[\gamma] + R_2[\gamma,\gamma]+\dots$ whose coefficients are homogeneous polynomial maps from $\Gamma$ into vector fields on $M$ such that
\begin{itemize}
\item $R_0$ is the limiting roto-rate.
\item $F_\gamma^*R_\gamma = R_\gamma$ in the sense of formal power series.
\item $\exp(2\pi\mathcal{L}_{R_\gamma})=1$ in the sense of formal power series.
\end{itemize}
\end{definition}

%\todo{Why is $G_\gamma$ a $O(\gamma)$ formal power series as opposed to just a formal commutator series in $\gamma$? In particular, in the definition, it is not approximating anything yet.}

\begin{definition}
Let $G_\gamma =G_1[\gamma]+ G_2[\gamma,\gamma] + \dots $ be an $O(\gamma)$ (no constant term) formal power series whose coefficients are homogeneous polynomial maps from $\Gamma$ into vector fields on $M$.
The \textbf{Lie transform} with \textbf{generator} $G_\gamma$ is the formal power series $\exp(\mathcal{L}_{G_\gamma})$ whose coefficients are homogeneous polynomial maps from $\Gamma$ into differential operators on the tensor algebra over $M$. 
\end{definition}

\begin{definition}\label{normalizing_transformation_map_def}
A \textbf{normalizing transformation} for a nearly-periodic map $F$ with roto-rate $R_\gamma$ is a Lie transform $\exp(\mathcal{L}_{G_\gamma})$ with generator $G_\gamma$ such that $R_\gamma = \exp(\mathcal{L}_{G_\gamma})R_0$.
\end{definition}

Our first and most fundamental result concerning nearly-periodic maps establishes the existence and uniqueness of the roto-rate in the non-resonant case. Like the corresponding result due to Kruskal in continuous time, this result holds to all orders in perturbation theory. 

\begin{theorem}[Existence and uniqueness of the roto-rate]\label{existence_of_rr}
Each \underline{non-resonant} nearly-periodic map admits a unique roto-rate. 
\end{theorem}
\begin{proof}
First we will show that there exists a Lie transform with generator $G_\gamma$ such that $R_\gamma\equiv \exp(\mathcal{L}_{G_\gamma})R_0$ is a roto-rate. 

To that end, we introduce a convenient way of representing $\gamma$-dependent pullback operators at the level of formal power series. Let $\psi_\gamma$ be a smooth $\gamma$-dependent diffeomorphism on $M$. By the Lie derivative formula, there is a unique $\gamma$-dependent $\Gamma^*$-valued vector field $W_\gamma$ such that
\begin{align}
    \left.\frac{d}{ds}\right|_{s=0}\psi_{\gamma+s\,\delta \gamma}^* = \psi_\gamma^*\mathcal{L}_{\langle W_\gamma,\delta \gamma \rangle},\label{W_def}
\end{align}
for each $\gamma,\delta\gamma \in \Gamma$. Here $\langle \cdot,\cdot\rangle$ denotes the natural pairing between $\Gamma$ and its dual space $\Gamma^*$. The object $W_\gamma$ both determines and is determined by the pullback operator $\psi_\gamma^*$ at the level of formal power series. This follows from recursive application of the identity
\begin{align}
    \psi_{s\gamma}^* = \psi_0^* + \int_0^s\psi_{s_1\gamma}^*\mathcal{L}_{\langle W_{s_1\,\gamma},\gamma\rangle}\,ds_1,\label{recusion_formula}
\end{align}
which may be understood as a consequence of \eqref{W_def} and the fundamental theorem of calculus. This can be viewed as Picard iteration of~\eqref{W_def} or fixed point iteration of~\eqref{recusion_formula}. The first step in the recursion is to substitute \eqref{recusion_formula} with $s=s_1$ into \eqref{recusion_formula} with $s=1$, resulting in
\begin{align*}
    \psi_{\gamma}^* &= \psi_0 + \int_0^1\left(  \psi_0^* + \int_0^{s_1}\psi_{s_2\gamma}^*\mathcal{L}_{\langle W_{s_2\,\gamma},\gamma\rangle}\,ds_2  \right)\mathcal{L}_{\langle W_{s_1\,\gamma},\gamma\rangle}\,ds_1\\
    & = \psi_0 + \int_0^1\psi_0^*   \mathcal{L}_{\langle W_{s_1\,\gamma},\gamma\rangle}\,ds_1 
    + \int_0^1\int_0^{s_1}\psi_{s_2\gamma}^*\mathcal{L}_{\langle W_{s_2\,\gamma},\gamma\rangle} \mathcal{L}_{\langle W_{s_1\,\gamma},\gamma\rangle}\,ds_2\,ds_1.
\end{align*}
The second step involves substituting \eqref{recusion_formula} with $s=s_2$ into the preceding expression, thereby producing a triple integral.
Continuing in this manner, it is straightforward to derive the following time-ordered exponential formulas for both the pullback $\psi_\gamma^*$ and pushforward operator $\psi_{\gamma *}$,
\begin{align}
    \psi_\gamma^* & = \psi_0^*\left[1 + \int_0^1 \mathcal{L}_{\langle W_{s_1\,\gamma},\gamma\rangle}\,ds_1 + \int_0^1\int_0^{s_1}\mathcal{L}_{\langle W_{s_2\,\gamma},\gamma\rangle}\mathcal{L}_{\langle W_{s_1\,\gamma},\gamma\rangle}\,ds_2\,ds_1 +\dots \right]\label{pullback_time_ordered_exp}\\
    \psi_{\gamma *} & = \left[1 - \int_0^1 \mathcal{L}_{\langle W_{s_1\,\gamma},\gamma\rangle}\,ds_1 + \int_0^1\int_0^{s_1}\mathcal{L}_{\langle W_{s_1\,\gamma},\gamma\rangle}\mathcal{L}_{\langle W_{s_2\,\gamma},\gamma\rangle}\,ds_2\,ds_1+\dots\right]\,\psi_{0\,*}.\label{pushforward_time_ordered_exp}
\end{align}
Upon introducing the formal power series expansion $W_\gamma = W_0 + W_1[\gamma] + W_2[\gamma,\gamma] + \dots$, the integrals in these formulas can be carried out, leading to the somewhat more explicit formulas
\begin{align}
    \psi_\gamma^* & =\psi_0^*\left[1 + \mathcal{L}_{\langle W_0,\gamma \rangle} + \frac{1}{2}\mathcal{L}_{\langle W_1[\gamma],\gamma \rangle} + \frac{1}{2}\mathcal{L}^2_{\langle W_0,\gamma \rangle}+O(\gamma^3)\right]\label{pullback_explicit}\\
    \psi_{\gamma *} & = \left[1 - \mathcal{L}_{\langle W_0,\gamma \rangle} - \frac{1}{2}\mathcal{L}_{\langle W_1[\gamma],\gamma \rangle} + \frac{1}{2}\mathcal{L}^2_{\langle W_0,\gamma \rangle}+O(\gamma^3)\right]\,\psi_{0\,*}.\label{pushforward_explicit}
\end{align}

% For example, since
% \begin{align*}
%     \psi_{\gamma}^* &= \psi_0^* + \int_0^1\psi_{s_1\gamma}^*\mathcal{L}_{\langle W_{s_1\,\gamma},\gamma\rangle}\,ds_1\\
%     & = \psi_0^* + \int_0^1\left(\psi_0^* + \int_0^1\psi_{s_2\,s_1\,\gamma}^*\mathcal{L}_{\langle W_{s_2\,s_1\,\gamma},s_1\,\gamma\rangle}\,ds_2\right)\mathcal{L}_{\langle W_{s_1\,\gamma},\gamma\rangle}\,ds_1\\
%     & = \psi_0^* + \int_0^1 \psi_0^*\mathcal{L}_{\langle W_{s_1\,\gamma},\gamma\rangle }\,ds_1 + \int_0^1\int_0^{s_1}\psi_{s_2\gamma}^*\mathcal{L}_{\langle W_{s_2\,\gamma},\gamma\rangle}\mathcal{L}_{\langle W_{s_1\,\gamma},\gamma\rangle}\,ds_2\,ds_1,
% \end{align*}
% we have
% \begin{align*}
%     \psi_{\gamma}^* &= \psi_0^* + \psi_0^*\mathcal{L}_{\langle W_{0},\gamma\rangle} \\
%     & + \frac{1}{2}\psi_0^*\mathcal{L}_{\langle W_{1}[\gamma],\gamma \rangle} + \frac{1}{2}\psi_0^*\mathcal{L}_{\langle W_0,\gamma \rangle}^2+ O(\gamma^3),
% \end{align*}
% where we have introduced the coefficients in the formal power series expansion $W_\gamma = W_0 + W_1[\gamma] + W_2[\gamma,\gamma] + \dots$. More generally, each term in the power series expansion of $\psi_{\gamma}^*$ may be expressed in terms of the $W_k$.

The preceding discussion applies in particular to $\psi_\gamma^* = F_\gamma^*$. In this case we will use the symbol $V_\gamma$ for $W_\gamma$. The discussion also applies to the formal pullback operator $\psi_\gamma^* = \phi_\gamma^*$, where
\[
\phi_\gamma^* = \exp(\mathcal{L}_{G_\gamma}),
\]
as well as its inverse $\phi_{\gamma *} = (\phi_\gamma^*)^{-1}$.
In this case we will use $\xi_\gamma$ in place of $W_\gamma$. Thus, we have the defining identities
\begin{align}
    \frac{d}{ds}\bigg|_{s=0}F_{\gamma+s\,\delta \gamma}^* &= F_\gamma^*\mathcal{L}_{\langle V_\gamma,\delta \gamma\rangle}\\
    \frac{d}{ds}\bigg|_{s=0}\phi_{\gamma+s\,\delta \gamma}^* &= \phi_\gamma^*\mathcal{L}_{\langle \xi_\gamma,\delta \gamma\rangle}.\label{Lie_transform_reparam}
\end{align}

We will now establish existence of the Lie transform with generator $G_\gamma$ by constructing an appropriate $\xi_\gamma$. The Lie transform itself can then be constructed  using the formulas \eqref{pullback_time_ordered_exp} and \eqref{pushforward_time_ordered_exp}. Define $R_\gamma = \exp(\mathcal{L}_{G_\gamma})R_0 = \phi_\gamma^*R_0$, where $G_\gamma$, or equivalently $\xi_\gamma$, is yet to be determined. This $R_\gamma$ satisfies $\exp(2\pi\mathcal{L}_{R_\gamma})=1$ and $R_{\gamma=0}=R_0$ automatically. We will determine the formal power series $\xi_\gamma = \xi_0 + \xi_1[\gamma] +\xi_2[\gamma,\gamma]+\dots$ by requiring $F_\gamma^*R_\gamma = R_\gamma$. If this can be done then $R_\gamma$ will be a roto-rate. 

The equation we would like to solve is equivalent to
\begin{align}
    (\phi_\gamma^{-1})^*F_\gamma^*\phi_\gamma^*R_0 = R_0,\label{transformed_objective_eqn}
\end{align}
where $(\phi_\gamma^{-1})^* = (\phi_\gamma^*)^{-1}$. Formally, this is just the statement that the ``diffeomorphism" $\overline{F}_\gamma = \phi_\gamma\circ F_\gamma\circ \phi_\gamma^{-1}$ preserves the limiting roto-rate $R_0$. Instead of solving  \eqref{transformed_objective_eqn} directly, we will demand that its $\gamma$-derivative vanishes. This derivative condition is
\begin{align}
   0 =  \frac{d}{ds}\bigg|_{s=0}\overline{F}_{\gamma+s\,\delta \gamma}^*R_0 = \overline{F}_\gamma^*\mathcal{L}_{\langle\overline{V}_\gamma,\delta\gamma\rangle}R_0,
\end{align}
where $\overline{V}_\gamma$ is readily shown to be given by
\begin{align*}
    \overline{V}_\gamma = \xi_\gamma - \overline{F}_{\gamma\,*}\xi_\gamma + \phi_{\gamma\,*}V_\gamma.
\end{align*}
Note that requiring the $\gamma$-derivative of \eqref{transformed_objective_eqn} to vanish implies \eqref{transformed_objective_eqn} itself since the latter is clearly satisfied when $\gamma=0$. Also note that since $\overline{F}_\gamma^*$ is formally invertible, the derivative condition is equivalent to 
\begin{align}
    \mathcal{L}_{R_0}\overline{V}_\gamma = 0.\label{derivative_condition}
\end{align}

To solve \eqref{derivative_condition}, we will expand the equation in powers of $\gamma$ and then argue inductively that each equation in the resulting sequence can be solved. At $O(\gamma^0)$ we have
\begin{align*}
    \mathcal{L}_{R_0}(\xi_0 - \Phi_{\theta_0\,*}\xi_0 + V_0) = 0.
\end{align*}
Denoting the limiting $U(1)$-average operation as $\langle \bm{T} \rangle = (2\pi)^{-1}\int_0^{2\pi}\Phi_\theta^*\bm{T}\,d\theta$, where $\bm{T}$ is any tensor field on $M$, and $\bm{T}^{\text{osc}} = \bm{T} - \langle \bm{T}\rangle$, this equation is equivalent to
\begin{align*}
    A_{\theta_0}\xi_0^{\text{osc}} = -V_0^{\text{osc}},
\end{align*}
where we have introduced the \textbf{homological operator}
\begin{align*}
    A_{\theta_0} = 1 - \Phi_{\theta_0\,*}.
\end{align*}
Since $F_\gamma$ is assumed to be non-resonant, the homological operator, regarded as a linear automorphism of the oscillating subspace of vector fields, is invertible. We may therefore solve the $O(\gamma^0)$ equation by setting
\begin{align*}
    \xi_0 & = -A_{\theta_0}^{-1}V_0^{\text{osc}}.
\end{align*}
At $O(\gamma^n)$, \eqref{derivative_condition} leads to 
\begin{align*}
    A_{\theta_0}\xi_n[\gamma,\dots,\gamma]^{\text{osc}} = S_n[\gamma,\dots,\gamma]^{\text{osc}},
\end{align*}
where $S_n[\gamma,\dots,\gamma]$ involves coefficients of $V_\gamma = V_0 + V_1[\gamma] + V_2[\gamma,\gamma]$ and $\xi_\gamma = \xi_0 + \xi_1[\gamma] + \xi_2[\gamma,\gamma] + \dots$ with polynomial degree (for $\xi$) at most $n-1$. Assuming the $\xi_k$ with $k < n$ have already been determined by solving the $O(\gamma^k)$ components of \eqref{derivative_condition}, we may therefore solve the $O(\gamma^n)$ equation by setting
\begin{align*}
    \xi_n[\gamma,\dots,\gamma] = A_{\theta_0}^{-1}S_n[\gamma,\dots,\gamma]^{\text{osc}}.
\end{align*}
Since we have already established that the $O(\gamma^0)$ equation has a solution, we now conclude by induction that \eqref{derivative_condition} may be solved for $\xi_\gamma$ to all orders in $\gamma$. It follows that a roto-rate exists.

Next we prove uniqueness of the $R_\gamma$ just constructed. Suppose $R_\gamma^\prime$ is a possibly distinct roto-rate, and consider the commutator $C_\gamma = [R_\gamma,R_\gamma^\prime]$. Immediate properties of $C_\gamma$ include $C_0 = 0$ and $F_\gamma^*C_\gamma = C_\gamma$. By construction of $R_\gamma$, we have $R_\gamma = \exp(\mathcal{L}_{G_\gamma})R_0$, which implies $\overline{C}_\gamma = \exp(-\mathcal{L}_{G_\gamma})C_\gamma = [R_0,\overline{R}_\gamma^\prime]$, where $\overline{R}_\gamma^\prime = \exp(-\mathcal{L}_{G_\gamma})R_\gamma^\prime$. Thus, the mean $\langle \overline{C}_\gamma\rangle =0$ vanishes to all orders in $\gamma$. Since $F_\gamma^*C_\gamma = C_\gamma$, we also have $\overline{F}_\gamma^* \overline{C}_\gamma = \overline{C}_\gamma$, where $\overline{F}_\gamma^*\equiv \exp(-\mathcal{L}_{G_\gamma})F_\gamma^*\exp(\mathcal{L}_{G_\gamma}) $. The first-order consequence of the last condition is $\Phi_{\theta_0}^*\overline{C}_1[\gamma] = \overline{C}_1[\gamma]$, which implies $\mathcal{L}_{R_0}\overline{C}_1[\gamma] = 0$ by non-resonance. But since the mean of $\overline{C}_1[\gamma] $ vanishes, we must have $\overline{C}_1[\gamma]=0$ for all $\gamma\in \Gamma$. The same argument applied repeatedly implies $\overline{C}_n=0$ for all $n\geq 1$. In other words, $R_\gamma$ and $R_\gamma^\prime$ commute. To finish the argument for unniqueness, we now use commutativity of $R_\gamma$ and $R_\gamma^\prime$ to find $1=\exp(2\pi\mathcal{L}_{R_\gamma}) = \exp(2\pi\mathcal{L}_{R_\gamma^\prime} + 2\pi\mathcal{L}_{R_\gamma - R_\gamma^\prime}) = \exp(2\pi\mathcal{L}_{R_\gamma - R_\gamma^\prime}),$ which can only be satisfied if $R_\gamma - R_\gamma^\prime = 0$, since $R_0 - R_0^\prime=0$.

\end{proof}

\begin{theorem}[Existence of normalizing transformations]\label{normalizing_transformation_existence}
Each \underline{non-resonant} nearly-periodic map admits a normalizing transformation.
\end{theorem}
\begin{proof}
This follows immediately from the proof of Theorem \ref{existence_of_rr}.
\end{proof}

\subsection{Nearly-periodic maps with Hamiltonian structure}
As in the continuous-time theory, existence of the roto-rate leads to additional insights for nearly-periodic maps that are Hamiltonian in an appropriate sense. In this subsection, we will establish the basic properties of nearly-periodic maps with Hamiltonian structure. We start by defining what we mean by Hamiltonian structure.
\begin{definition}
A \textbf{$\gamma$-dependent presymplectic manifold} is a manifold $M$ equipped with a smooth $\gamma$-dependent $2$-form $\Omega_\gamma$ such that $\mathbf{d}\Omega_\gamma = 0$ for each $\gamma\in \Gamma$. We say $(M,\Omega_\gamma)$ is \textbf{exact} when there is a smooth $\gamma$-dependent $1$-form $\vartheta_\gamma$ such that $\Omega_\gamma = -\mathbf{d}\vartheta_\gamma$.
\end{definition}

\begin{definition}[Presymplectic nearly-periodic map]
A \textbf{Presymplectic nearly-periodic map} on a $\gamma$-dependent presymplectic manifold $(M,\Omega_\gamma)$ is a nearly-periodic map $F$ such that $F_\gamma^*\Omega_\gamma = \Omega_\gamma$ for each $\gamma\in \Gamma$.
\end{definition}

\begin{definition}[Hamiltonian nearly-periodic map]
A \textbf{Hamiltonian nearly-periodic map} on a $\gamma$-dependent presymplectic manifold $(M,\Omega_\gamma)$ is a nearly-periodic map $F$ such that there is a smooth $(t,\gamma)$-dependent vector field $Y_{t,\gamma}$ with the following properties: 
\begin{itemize}
\item $t\in \mathbb{R}$.
\item $\iota_{Y_{t,\gamma}}\Omega_\gamma = \mathbf{d}H_{t,\gamma}$, for some smooth $(t,\gamma)$-dependent function $H_{t,\gamma}$.
\item For each $\gamma\in\Gamma$, $F_\gamma$ is the $t=1$ flow of $Y_{t,\gamma}$.
\end{itemize} 

\end{definition}

\begin{lemma}
Each Hamiltonian nearly-periodic map is a presymplectic nearly-periodic map.
\end{lemma}

\begin{theorem}[Roto-rate is presymplectic]\label{rr_is_presymplectic}
If $F$ is a \underline{non-resonant} presymplectic nearly-periodic map on a $\gamma$-dependent presymplectic manifold $(M,\Omega_\gamma)$ with roto-rate $R_\gamma$ then $\mathcal{L}_{R_\gamma}\Omega_\gamma = 0$.
\end{theorem}
\begin{proof}
First note that presymplecticity of $F$ with $\gamma=0$ implies $F_0^*\Omega_0 =\Phi_{\theta_0}^*\Omega_0= \Omega_0$. Upon introducing the $2$-form-valued Fourier coefficients,
\begin{align}
\Omega_0^k = \frac{1}{2\pi}\int_0^{2\pi}e^{-ik\theta}\Phi_\theta^*\Omega_0\,d\theta,\quad k\in\mathbb{Z},\label{Fourier_coefficients}
\end{align}
the last identity may be rewritten as the sequence of identities $e^{ik\theta_0}\Omega_0^k = \Omega_0^k$, $k\in\mathbb{Z}$. But by non-resonance of $F$, $1-e^{ik\theta_0}$ is non-vanishing for each $k$. We conclude that $\Omega_0^k=0$ for non-zero $k$, or, equivalently, $\mathcal{L}_{R_0}\Omega_0 = 0$.

Presymplecticity of $F$ for non-zero $\gamma$ implies $F_\gamma^*\Omega_\gamma = \Omega_\gamma$ for each $\gamma\in \Gamma$. Applying the Lie derivative $\mathcal{L}_{R_\gamma}$ to this identity and using $F_\gamma^*R_\gamma = R_\gamma$ implies $F_\gamma^*(\mathcal{L}_{R_\gamma}\Omega_\gamma) =(\mathcal{L}_{R_\gamma}\Omega_\gamma) $. In other words, $\alpha_\gamma = \mathcal{L}_{R_\gamma}\Omega_\gamma$ is (formally) an absolute integral invariant for $F_\gamma$. By the argument from the previous paragraph, we see immediately that $\alpha_0 = 0$. To finish the proof, we will use integral invariance together with existence of a normalizing transformation to find that $\alpha_\gamma=0$ to all orders in $\gamma$. This argument will parallel the proof of Proposition \ref{bootstrap_prop}.

Let $G_\gamma$ be the generator of a normalizing transformation for $F$ given by Theorem \ref{normalizing_transformation_existence}. Set $\overline{\alpha}_\gamma = \exp(-\mathcal{L}_{G_\gamma})\alpha_\gamma = \overline{\alpha}_0 + \overline{\alpha}_1[\gamma] +\overline{\alpha}_2[\gamma,\gamma] +\dots$. Note that $\alpha_\gamma = 0$ if and only if $\overline{\alpha}_\gamma=0$. Since $\alpha_\gamma$ is an integral invariant for $F_\gamma$, $\overline{\alpha}_\gamma$ must satisfy \begin{align}\exp(-\mathcal{L}_{G_\gamma})F_\gamma^*\exp(\mathcal{L}_{G_\gamma})\overline{\alpha}_\gamma = \overline{\alpha}_\gamma.\label{transformed_iv}\end{align} Because $\overline{\alpha}_0 = \alpha_0 = 0$, the first-order consequence of \eqref{transformed_iv} is $F_0^*\overline{\alpha}_1[\gamma] = \Phi_{\theta_0}^*\overline{\alpha}_1[\gamma] =\overline{\alpha}_1[\gamma] $. By our earlier argument, we must then have $(\overline{\alpha}_1[\gamma])^k=0$ for $k\neq 0$. But using $\exp(-\mathcal{L}_{G_\gamma})R_\gamma = R_0$, we also find $\overline{\alpha}_\gamma = \mathcal{L}_{R_0}\overline{\Omega}_\gamma$, where $\overline{\Omega}_\gamma =  \exp(-\mathcal{L}_{G_\gamma})\Omega_\gamma$. The latter implies $\overline{\alpha}_\gamma^0 = 0$, and $(\overline{\alpha}_1[\gamma])^0=0$ in particular. Thus, $\overline{\alpha}_1[\gamma]=0$ for all $\gamma\in\Gamma$. We may now repeat this argument for $\overline{\alpha}_2[\gamma,\gamma]$, $\overline{\alpha}_3[\gamma,\gamma,\gamma]$, etc., to obtain the desired result.

\end{proof}

Using presymplecticity of the roto-rate, we may now use a version of Noether's theorem to establish existence of adiabatic invariants for many interesting presymplectic nearly-periodic maps.

\begin{theorem}[Existence of an adiabatic invariant]\label{nps_noether}
Let $F$ be a \underline{non-resonant} presymplectic nearly-periodic map on the exact $\gamma$-dependent presymplectic manifold $(M,\Omega_\gamma)$ with roto-rate $R_\gamma$. Assume one of the following conditions is satisfied.
\begin{itemize}
\item[(1)] $F$ is Hamiltonian.
\item[(2)] $M$ is connected and the limiting roto rate $R_0$ has at least one zero.
\end{itemize} 
There exists a smooth $\gamma$-dependent $1$-form $\theta_\gamma$ such that $\mathcal{L}_{R_\gamma}\theta_\gamma = 0$ and $-\mathbf{d}\theta_\gamma =\Omega_\gamma$ in the sense of formal power series. Moreover the quantity
\begin{align}
\mu_\zeta = \iota_{R_\gamma}\theta_\gamma\label{the_adiabatic_invariant}
\end{align}
satisfies $F_\gamma^*\mu_\gamma = \mu_\gamma$ in the sense of formal power series. In other words, $\mu_\gamma$ is an adiabatic invariant for $F$.
\end{theorem}
\begin{proof}
By Theorem \ref{rr_is_presymplectic}, a primitive $\theta_\gamma$ with the desired properties may be constructed as in the proof of Proposition \ref{existence_of_mu}.

To establish adiabatic invariance of $\mu_\gamma$, first we compute the exterior derivative of $\mu_\gamma$ using Cartan's formula to obtain $\mathbf{d}\mu_\gamma = \iota_{R_\gamma}\Omega_\gamma$. Since both $R_\gamma$ and $\Omega_\gamma$ are $F_\gamma$-invariant, it follows that $\mathbf{d}F_\gamma^*\mu_\gamma = \mathbf{d}\mu_\gamma$. The difference $c_\gamma \equiv F_\gamma^*\mu_\gamma - \mu_\gamma$ must therefore be a formal power series whose coefficients are homogeneous polynomial maps from $\Gamma$ into locally-constant functions on $M$. To complete the proof, we must demonstrate that $c_\gamma = 0$ to all orders.

First suppose $F$ is Hamiltonian. Then there is a smooth $(t,\gamma)$-dependent Hamiltonian vector field $Y_{t,\gamma}$ with Hamiltonian $H_{t,\gamma}$ whose $t=1$ flow is equal to $F_\gamma$. Let $\mathcal{F}^\gamma_t$ denote the time-$t$ flow map for $Y_{t,\gamma}$ with $\mathcal{F}^\gamma_0 = \text{id}_M$. By the fundamental theorem of calculus, the definition of Lie derivative, and Cartan's formula, we therefore have
\begin{align*}
F_\gamma^*\mu_\gamma &= \iota_{R_\gamma}(\mathcal{F}^\gamma_1)^*\theta_\gamma\\
& =  \iota_{R_\gamma}\theta_\gamma +\iota_{R_\gamma}\int_0^1 \frac{d}{dt} (\mathcal{F}^\gamma_t)^*\theta_\gamma\,dt\\
& = \mu_\gamma + \iota_{R_\gamma}\int_0^1 (\mathcal{F}^\gamma_t)^*(\mathcal{L}_{Y_{t,\gamma}}\theta_\gamma)\,dt\\
& = \mu_\gamma + \iota_{R_\gamma}\int_0^1 (\mathcal{F}^\gamma_t)^*(\iota_{Y_{t,\gamma}}\mathbf{d}\theta_\gamma + \mathbf{d}\iota_{Y_{t,\gamma}}\theta_\gamma )\,dt\\
& = \mu_\gamma + \iota_{R_\gamma}\int_0^1 (\mathcal{F}^\gamma_t)^* \mathbf{d}(\iota_{Y_{t,\gamma}}\theta_\gamma -H_{t,\gamma})\,dt\\
& = \mu_\gamma + \mathcal{L}_{R_\gamma}\int_0^1 (\mathcal{F}^\gamma_t)^* (\iota_{Y_{t,\gamma}}\theta_\gamma -H_{t,\gamma})\,dt.
\end{align*}
Applying $\exp(\theta\mathcal{L}_{R_0})$ to this identity and averaging in $\theta$ gives the desired result.

Finally suppose that $M$ is connected and that $R_0(m_0)=0$ for some $m_0\in M$.  Let $G_\gamma$ be the generator of a normalizing transformation. We have
\begin{align*}
\exp(-\mathcal{L}_{G_\gamma})\mu_\gamma(m_0) = \iota_{R_0}\exp(-\mathcal{L}_{G_\gamma}){\theta}_\gamma(m_0)=0,
\end{align*}
and
\begin{align*}
\exp(-\mathcal{L}_{G_\gamma})F_\gamma^*\mu_\gamma(m_0) &= \exp(-\mathcal{L}_{G_\gamma})\iota_{R_\gamma}F_\gamma^*\theta_\gamma(m_0) = \iota_{R_0}(\exp(-\mathcal{L}_{G_\gamma})F_\gamma^*\theta_\gamma)(m_0)=0.
\end{align*}
It follows that $c_\gamma$ is zero on the connected component of $M$ containing $m_0$. But since $M$ is connected, $c_\gamma$ is therefore zero everywhere, as claimed.

\end{proof}

\begin{remark}
A simple example illustrates how existence of an adiabatic invariant can fail. Let $M =S^1\times\mathbb{R} \ni(\zeta,I)$, $\Omega_\gamma = d\zeta\wedge dI = -\mathbf{d}(I\,d\zeta)$, and $\Gamma =\mathbb{R}$. The mapping $F(\zeta,I,\gamma) = (\zeta+\theta_0,I+\gamma)$ defines a non-resonant nearly-periodic map for almost all $\theta_0\in U(1)$. The roto-rate is given to all orders by $R_\gamma = \partial_\zeta$. Moreover, $F_\gamma$ is area-preserving for each $\gamma$, and hence presymplectic. The quantity $\mu_\gamma = \iota_{R_\gamma}(I\,d\zeta) = I$ from \eqref{the_adiabatic_invariant} is clearly not an adiabatic invariant for $F$ since $F_\gamma^*I = I + \gamma$.  Note that, in this example, the $R_0$-orbits are not contractible and that $F$ is presymplectic but not Hamiltonian.
\end{remark}

\subsection{Geometric integration of noncanonical Hamiltonian systems using nearly-periodic maps}\label{sec:symplectic_lorentz}

Let $Q\subset E$ be a connected open subset of a finite-dimensional vector space $E$ equipped with an exact symplectic form $\omega = -\mathbf{d}\vartheta$. Consider a Hamiltonian system on $(Q,\omega=-\mathbf{d}\vartheta)$ with Hamiltonian $H:Q\rightarrow\mathbb{R}$. Without loss of generality \cite{Silva_book_2008,Burby_Hirvijoki_2021_JMP}, assume that $Q$ is equipped with an almost complex structure $\mathbb{J}:TQ\rightarrow TQ$ compatible with $\omega = -\mathbf{d}\vartheta$, so that $g(v,w) = \omega(v,\mathbb{J}w)$ defines a metric tensor on $Q$. In this scenario, we may equip the tangent bundle $\pi:TQ\rightarrow Q$ with the ``magnetic" symplectic form
\begin{align*}
\Omega_\epsilon^* = \pi^*\omega +\epsilon\,\Omega,
\end{align*}
where $\epsilon$ is a real parameter and $\Omega$ is the pullback of the canonical symplectic form on $T^*Q$ along the bundle map $TQ\rightarrow T^*Q$ defined by $g$. We may also define a natural Hamilton function on $TQ$,
\begin{align*}
H^*_\epsilon(q,v) =\frac{1}{2} \epsilon^2\,g_q(v,v) + \epsilon\,H(q).
\end{align*}
As explained in detail in \cite{Burby_Hirvijoki_2021_JMP}, $H^*_\epsilon$ defines a Hamiltonian nearly-periodic system whose slow manifold dynamics recovers the dynamics of $H$ on $Q$ as $\epsilon\rightarrow 0$. The limiting roto-rate is 
\begin{align}
R_0(q,v) = \mathbb{J}_q(v-X_H(q))\cdot\partial_v,\label{target_roto_sl}
\end{align}
where $X_H$ denotes the Hamiltonian vector field on $Q$ associated with $H$, and the angular frequency function $\omega_0 = 1$. Moreover, the adiabatic invariant associated with $H^*_\epsilon$ ensures that the slow manifold enjoys long-term normal stability. It is crucial the metric $g$ is determined by a almost complex structure $\mathbb{J}$ compatible with $\omega$ for these results to hold. If $g$ is a more general metric tensor then the larger system on $TQ$ need not be nearly-periodic, and an adiabatic invariant need not exist.

The purpose of this section is to combine observations from \cite{Burby_Hirvijoki_2021_JMP} with the theory of nearly-periodic maps in order to construct a geometric numerical integrator for $H$. The integrator will be given as an implicitly-defined mapping on $TQ$ that is provably presymplectic nearly-periodic with limiting roto-rate $R_0$. We will show that this mapping admits a slow manifold diffeomorphic to $Q$ on which iterations of the map approximate the $H$-flow. In fact, the mapping is a slow manifold integrator in the sense described in \cite{Burby_Klotz_2020}. In addition, we will argue using the Noether theorem for nearly-periodic maps that this discrete-time slow manifold enjoys long-term normal stability. This ensures that the mapping on $TQ$ will function effectively and reliably as a structure-preserving integrator for the original Hamiltonian system on $Q$. We remark that the results described in this section provide a general solution to the problem of structure-preserving integration of non-canonical Hamiltonian systems on exact symplectic manifolds. For a completely different approach that is less geometric, we refer readers to \cite{Kraus_2017_arXiv}.

We begin with some preliminary remarks. 
\begin{remark}
It will be convenient to work with the parameter $\hbar = \sqrt{\epsilon}$ instead of $\epsilon$. There are technical reasons for doing so that will not be discussed here; however, an obvious physical benefit will be that $\hbar$ may be interpreted as a timestep. The symplectic form on $TQ$ will therefore be given by 
\begin{align*}
\Omega^*_{\hbar} = \pi^*\omega + \hbar^2\,\Omega.
\end{align*}
\end{remark}
\begin{remark}
It is useful to describe the goal of our construction in more concrete terms. We aim to find a smooth $\hbar$-dependent mapping $\Psi_{\hbar}:TQ\rightarrow TQ$ that is both non-resonant nearly-periodic with limiting roto-rate $R_0$ given by \eqref{target_roto_sl} and symplectic, $\Psi_{\hbar}^*\Omega^*_{\hbar} = \Omega^*_{\hbar}$, for all $\hbar\ll 1$. Since $Q$ is connected and $R_0$ has a manifold of fixed points of the form $\{(q,X_H(q))\}\subset TQ$, Theorem \ref{nps_noether} (Noether's theorem for nearly-periodic maps) ensures that the mapping we seek will admit an adiabatic invariant $\mu_{\hbar}$.  
\end{remark}

\begin{remark}
We may determine the leading-order term in the formal power series $\mu_{\hbar} = \mu_0 + \hbar\,\mu_1 +\hbar^2\,\mu_2 + \dots$ using only the explicit expressions for $\Omega^*_{\hbar}$ and $R_0$ in conjunction with the general existence theorem (Theorem \ref{nps_noether}). Recall that the theorem says that the adiabatic invariant is given by $\mu_{\hbar} = \iota_{R_{\hbar}}\overline{\Theta}_{\hbar}$, where $R_{\hbar}$ is the roto rate and $\overline{\Theta}_{\hbar}$ is a $U(1)$-invariant primitive for $\Omega^*_{\hbar}$. In particular, the roto-rate must satisfy Hamilton's equation $\iota_{R_{\hbar}}\Omega^*_{\hbar} = \mathbf{d}\mu_{\hbar}$ with Hamiltonian $\mu_{\hbar}$. We apply this theorem as follows.

If $f:TQ\rightarrow \mathbb{R}$ is any smooth $\hbar$-independent function on $TQ$ then it is straightforward to verify that $X_f=\hbar^{-4}(\mathbb{J}\,\partial_{v}f)\cdot\partial_v+\text{l.o.t.}$. Therefore we must have $R_{\hbar}=\hbar^{-4}(\mathbb{J}\,\partial_{v}\mu_0)\cdot\partial_v+\text{l.o.t.}$. But since $R_{\hbar}=O(1)$, the last expression implies $\partial_v\mu_0=0$ everywhere on $TQ$, or $\mu_0=\mu_0(q)$. In fact $\mu_0(q)$ must be constant, for if $(q,v_{\hbar}(q))$ is a zero for $R_{\hbar}$ then Hamilton's equation implies $0=\lim_{\hbar\rightarrow 0}\mathbf{d}\mu_{\hbar}(q,v_{\hbar}(q)) =  \mathbf{d}\mu_0(q)$. And by evaluating the formula $\mu_{\hbar} = \iota_{R_{\hbar}}\overline{\Theta}_{\hbar}$ at $(q,v_{\hbar}(q))$ we find that the constant must in fact be $0$. In other words $\mu_0 = 0$. Essentially the same argument may now be applied to conclude $\mu_1=\mu_2=\mu_3=0$; the argument for $\mu_1$ proceeds as follows. Since $\mu_0=0$, Hamilton's equation implies $R_0 = \hbar^{-3}(\mathbb{J}\,\partial_{v}\mu_1)\cdot\partial_v+\text{l.o.t.}$, or $\mu_1 = \mu_1(q)$. Evaluating Hamilton's equation at $(q,v_{\hbar}(q))$, dividing by $\hbar$, and taking the limit $\hbar\rightarrow 0$ then leads to $0 = \lim_{\hbar\rightarrow 0}\hbar^{-1}\mathbf{d}\mu_{\hbar}(q,v_{\hbar}(q)) =\mathbf{d}\mu_1(q)$, implying $\mu_1$ is a constant. Applying the same procedure to the formula $\mu_{\hbar} = \iota_{R_{\hbar}}\overline{\Theta}_{\hbar}$ implies the constant must be $0$.

We have now established that the adiabatic invariant for the nearly-periodic map we aim to construct must have the form $\mu_{\hbar} =\hbar^4\,\mu_4 + \hbar^5\,\mu_5 + \dots$. We can determine an explicit expression for $\mu_4$ as follows. By the above remarks, we must have $R_{\hbar}=(\mathbb{J}\,\partial_{v}\mu_4)\cdot\partial_v+\text{l.o.t.}$, which implies in particular that $R_0 = (\mathbb{J}\,\partial_{v}\mu_4)\cdot\partial_v$. Since the desired form of $R_0$ is known, we therefore obtain the following partial differential equation for $\mu_4$:
\begin{align*}
\mathbb{J}(v-X_H) = \mathbb{J}\,\partial_{v}\mu_4.
\end{align*}
The general solution of this equation is given by $\mu_4(q,v) = \frac{1}{2}g_q(v-X_H(q),v-X_H(q))+\chi(q)$, where $\chi(q)$ is an arbitrary function of $q$. To determine $\chi$, we evaluate the formula $\mu_{\hbar} = \iota_{R_{\hbar}}\overline{\Theta}_{\hbar}$ at a fixed point $(q,v_{\hbar}(q))$ to find $0 =\lim_{\hbar\rightarrow 0}\hbar^{-4}\mu_{\hbar}(q,v_{\hbar}(q)) = \mu_4(q,v_0(q))=\chi(q)$. Note that we have used the formula for $v_0(q) = X_H(q)$ for fixed points of $R_0$. We conclude that the adiabatic invariant must have the general form
\begin{align}
\mu_{\hbar}(q,v) =\hbar^{4}\frac{1}{2}g_q(v-X_H(q),v-X_H(q)) + O(\hbar^{5}). \label{mu_gen_prediction}
\end{align}
This formula will be useful later when we argue for long-term normal stability.
\end{remark}

To construct the mapping $\Psi_{\hbar}:TQ\rightarrow TQ$, we begin by introducing the generating function
\begin{align}
    S(q,\overline{q}) &= \int_q^{\overline{q}}\vartheta + \Phi^*\Sigma(q,\overline{q})  ,
\end{align}
where $\Phi(q,\overline{q}) = (q/2+\overline{q}/2,\overline{q}-q)$ is a diffeomorphism $Q\times Q\rightarrow TQ$, and $\Sigma:TQ\rightarrow \mathbb{R}$ is given by
\begin{align*}
\Sigma(x,\xi) &= -\hbar H(x) + \hbar^2\,\langle X_H(x),\xi \rangle  -\frac{1}{12}\hbar^2 \partial_k\omega_{j\ell}(x)\,X_H^k(x)\,X_H^j(x)\,\xi^\ell\nonumber\\
    & \qquad- \frac{1}{4}\left(\frac{\sin\theta_0}{1-\cos\theta_0}\right)\langle \xi-\hbar X_H(x),\xi-\hbar X_H(x) \rangle.
\end{align*}
Here, $\langle\cdot,\cdot\rangle$ is shorthand for $g(\cdot,\cdot)$, the integral is taken along the straight line connecting $q$ with $\overline{q}$, $X_H = -\mathbb{J}\nabla H$ is the Hamiltonian vector field associated with $H$, and $\theta_0\nin \{0,\pi\}$. The metric tensor, the Hamiltonian, and the Hamiltonian vector field are evaluated at the midpoint $x=(\overline{q}+q)/2$.

%\begin{align}
%    S(q,\overline{q}) &= -\hbar H + \hbar^2\,\langle X_H,\overline{q}-q \rangle  - \frac{1}{2}\hbar^3\,\langle X_H,X_H \rangle\nonumber\\
%    &+\int_q^{\overline{q}}\vartheta - \frac{1}{4}\left(\frac{\sin\theta_0}{1-\cos\theta_0}\right)\langle \overline{q}-q-\hbar X_H,\overline{q}-q-\hbar X_H \rangle,
%\end{align}
%
%\begin{align*}
%\Sigma(x,\xi) &= -\hbar H(x) + \hbar^2\,\langle X_H(x),\xi \rangle  - \frac{1}{2}\hbar^3\,\langle X_H(x),X_H(x) \rangle\nonumber\\
%    & - \frac{1}{4}\left(\frac{\sin\theta_0}{1-\cos\theta_0}\right)\langle \xi-\hbar X_H(x),\xi-\hbar X_H(x) \rangle,
%\end{align*}

\begin{definition}
The \textbf{symplectic Lorentz map} is the mapping $\Psi_{\hbar}:TQ\rightarrow TQ:(q,v)\mapsto (\overline{q},\overline{v})$ defined by the implicit relations
\begin{align}
    \vartheta_{\overline{q}} + \hbar^2\,g_{\overline{q}}(\overline{v},d\overline{q}) &= \mathbf{d}^{(1)}S\label{slm_fwd},\\
    \vartheta_q+ \hbar^2\,g_{{q}}({v},d{q}) &= -\mathbf{d}^{(2)}S.\label{slm_bkwd}
\end{align}
\end{definition}

\begin{proposition}\label{slm_existence}
The symplectic Lorentz map is well-defined and smooth in $(q,v,\hbar)$ for $\hbar$ in a neighborhood of $0\in\mathbb{R}$. Moreover, it preserves the $\hbar$-dependent symplectic form $\Omega^*_{\hbar}$ and satisfies 
\begin{align}
\Psi_0(q,v) = (x,X_H(q)+\exp(-\theta_0\,\mathbb{J}(q))[v-X_H(q)]) .\label{limit_symplectic_lorentz}
\end{align}
\end{proposition}
\begin{proof}
First we will construct a convenient moving frame on $Q\times Q$ onto which we will resolve the implicit relations \eqref{slm_fwd}-\eqref{slm_bkwd}. We will start by building a frame on $TQ$ and then finish by pulling back along the mapping $\Phi:Q\times Q\rightarrow TQ$ defined above. Without loss of generality, assume $Q=\mathbb{R}^n$ for an even integer $n$ and let $(x^i,\xi^i)$ denote the standard linear coordinate system on $TQ$. Fix $(x,\xi)\in TQ$ and let $\gamma:[-1,1]\rightarrow Q$ be a smooth curve in $Q$ with $\gamma(0) = x$. Relative to the Riemannian structure defined by the metric $g$ there is a unique horizontal lift $\tilde{\gamma}:[-1,1]\rightarrow TQ$ with $\tilde{\gamma}(0) = (x,\xi)$. In this manner, to each $(x,\xi)\in TQ$ and each tangent vector $w\in T_xQ$ we assign a lifted tangent vector $\tilde{w}\in T_{(x,\xi)}TQ$. Applying this construction point-wise to the coordinate vector fields $\partial_{x^i}$ on $Q$, we obtain linearly-independent vector fields $\tilde{\partial}_{x^i}$ on $TQ$. The collection of $2n$ vector fields $(\tilde{\partial}_{x^i},\partial_{\xi^i})$ comprise a frame on $TQ$. A frame on $Q\times Q$ is then furnished by the vector fields $(U^i,A^i)$, where $U^i = \Phi^*\tilde{\partial}_{x^i}$ and  $\quad A^i  =\Phi^*\partial_{\xi^i}$.
Upon introducing the Christoffel symbols $\nabla_{\partial_{x^i}}\partial_{x^j} = \Gamma^k_{\,\,ij}\,\partial_{x^{k}}$, the vector fields $\tilde{\partial}_{x^i}$ may be written as
\begin{align*}
\tilde{\partial}_{x^i}  = \partial_{x^i} -\Gamma^k_{\,\,ij}(x)\,\xi^j\,\partial_{\xi^k}.
\end{align*}
Therefore, we may write the following explicit formulas for the frame $(U^i,A^i)$:
\begin{align*}
U^i 
&=\bigg(\partial_{q^i} +\frac{1}{2}\Gamma^k_{\,\,ij}(x)\,\xi^j\,\partial_{q^k}\bigg)+ \bigg(\partial_{\overline{q}^i} -  \frac{1}{2}\Gamma^k_{\,\,ij}(x)\,\xi^j\,\partial_{\overline{q}^k}\bigg)\\
A^i & = -\frac{1}{2}\partial_{q^i} + \frac{1}{2}\partial_{\overline{q}^i},
\end{align*}
where we remind the reader that $x = (q+\overline{q})/2$ and $\xi = \overline{q}-q$.

%\begin{align*}
%\Phi^{-1}(x,\xi) &= (x-\xi/2,x+\xi/2)\\
%\Phi^*\partial_{x^i} & = \partial_{q^i} + \partial_{\overline{q}^i}\\
%\Phi^*\partial_{\xi^i} & = -\frac{1}{2}\partial_{q^i} + \frac{1}{2}\partial_{\overline{q}^i}
%\end{align*}

Next we re-write \eqref{slm_fwd}-\eqref{slm_bkwd} as a single equation on $Q\times Q$,
\begin{align}
\vartheta_{\overline{q}}(d\overline{q}) -\vartheta_{q}(dq)+ \hbar^2\,g_{\overline{q}}(\overline{v},d\overline{q}) -  \hbar^2\,g_{{q}}({v},d{q})= \mathbf{d}S,\label{proof_symp_condition}
\end{align}
and then take components along the frame $(U^i,A^i)$ to obtain 
\begin{align*}
&\hbar^2\,g_{\ell\,i}(\overline{q})\,\overline{v}^\ell - \hbar^2\,g_{\ell\,i}(q)\,v^\ell -\frac{1}{2} \hbar^2\,\Gamma^k_{\,\,ij}(x)\,\xi^j\bigg(g_{\ell\,k}(\overline{q})\,\overline{v}^\ell + g_{\ell\,k}(q)\,v^\ell\bigg)\\
&=\int_0^1\xi^j\omega_{ji}(\lambda) \,d\lambda-\int_0^1\left[\lambda-\frac{1}{2}\right] \xi^j\,\omega_{jk}(\lambda)\Gamma^k_{\,\,i\ell}(x)\xi^\ell\,d\lambda\nonumber\\
&\qquad-\hbar\,\partial_iH(x) +\frac{\hbar}{2}\left(\frac{\sin\theta_0}{1-\cos\theta_0}\right) \,g_{\ell\,j}(x)\,X^\ell_{H;i}(x)\,(\xi^j - \hbar\,X_H^j(x))\\
&\qquad-\frac{1}{12}\hbar^2\,\partial_i(\partial_k\omega_{j\ell}\,X_H^k\,X_H^j)\,\xi^\ell + \frac{1}{12}\hbar^2\,\partial_k\omega_{j\ell}\,X_H^k\,X_H^j\,\Gamma^\ell_{\,\,im}\,\xi^m,
\end{align*}
and
\begin{align*}
&\frac{1}{2}\,\hbar^2\,g_{\ell\,i}(\overline{q})\,\overline{v}^\ell + \frac{1}{2}\,\hbar^2\,g_{\ell\,i}(q)\,v^\ell\\
& = \int_0^1\left[\lambda-\frac{1}{2}\right]\,\xi^j\omega_{ji}(\lambda)\,d\lambda-\frac{1}{12}\hbar^2\,\partial_k\omega_{ji}\,X_H^k\,X_H^j \nonumber\\
& \qquad + \hbar^2\,g_{ij}(x)\,X_H^j(x) - \frac{1}{2}\left(\frac{\sin\theta_0}{1-\cos\theta_0}\right)\,g_{ij}(x)(\xi^j -\hbar\,X_H^j(x)).
\end{align*}
Here, we have applied the useful formulas
\begin{align*}
\mathbf{d}\int_{q}^{\overline{q}}\vartheta = \vartheta_{\overline{q}}(d\overline{q}) -\vartheta_q(dq)+\int_0^1\omega_\lambda(\xi,[1-\lambda]dq + \lambda\,d\overline{q})\,d\lambda
\end{align*}
and
\begin{align*}
\mathcal{L}_{\tilde{\partial}_{x^i}}\Sigma & = -\hbar\,\partial_iH +\frac{\hbar}{2}\left(\frac{\sin\theta_0}{1-\cos\theta_0}\right) \,g(\nabla_{\partial_i}X_H,\xi-\hbar\,X_H) \\
& \qquad-\frac{1}{12}\hbar^2\,\partial_i(\partial_k\omega_{j\ell}\,X_H^k\,X_H^j)\,\xi^\ell + \frac{1}{12}\hbar^2\,\partial_k\omega_{j\ell}\,X_H^k\,X_H^j\,\Gamma^\ell_{\,\,im}\,\xi^m\\
\mathcal{L}_{\partial_{\xi^i}}\Sigma & = \hbar^2\,g(X_H,\partial_{x^i})-\frac{1}{12}\hbar^2\,\partial_k\omega_{ji}\,X_H^k\,X_H^j - \frac{1}{2}\left(\frac{\sin\theta_0}{1-\cos\theta_0}\right)\,g(\xi-\hbar\,X_H,\partial_{x^i}),
\end{align*}
where $\omega_\lambda = \omega_{q(\lambda)}$ with $q(\lambda) = [1-\lambda]\,q + \lambda\,\overline{q}$.

To show that these implicit equations define a smooth $\hbar$-dependent mapping $\Psi_{\hbar}$, we first introduce the new variable $\Delta = \hbar^{-2}(\xi -\hbar\,X_H(x))$ and then observe that when expressed in terms of $\Delta$ the implicit equations above may be written in the form
\begin{align*}
&\hbar^2\,g_{\ell\,i}(x)\,(\overline{v}^\ell-v^\ell)=\hbar^2\Delta^j\omega_{ji}(x)+O(\hbar^3) 
\end{align*}
and
\begin{align*}
&\frac{1}{2}\,\hbar^2\,g_{\ell\,i}(x)\,(\overline{v}^\ell+v^\ell) =\hbar^2\,g_{ij}(x)\,X_H^j(x) - \hbar^2\frac{1}{2}\left(\frac{\sin\theta_0}{1-\cos\theta_0}\right)\,g_{ij}(x)\Delta^j+O(\hbar^3).
\end{align*}
Note in particular that the term $\frac{1}{12}\hbar^2\,\partial_k\omega_{ji}\,X_H^k\,X_H^j$ exactly cancels the second-order part of $ \int_0^1\left[\lambda-\frac{1}{2}\right]\,\xi^j\omega_{ji}(\lambda)\,d\lambda$. Dividing these expression by $\hbar^2$ implies that there are smooth functions $Z_1,Z_2:\mathbb{R}^n\times\mathbb{R}^n\times\mathbb{R}^n\times\mathbb{R}^n\times\mathbb{R}\rightarrow\mathbb{R}^n$ given by
\begin{align*}
Z_{1i}(x,v,\Delta,\overline{v},\hbar) &= g_{\ell\,i}(x)\,(\overline{v}^\ell-v^\ell) - \Delta^j\omega_{ji}(x) + O(\hbar),\\
Z_{2i}(x,v,\Delta,\overline{v},\hbar) & =\frac{1}{2}\,g_{\ell\,i}(x)\,(\overline{v}^\ell+v^\ell)  - \,g_{ij}(x)\,X_H^j(x) + \frac{1}{2}\left(\frac{\sin\theta_0}{1-\cos\theta_0}\right)\,g_{ij}(x)\Delta^j+O(\hbar),
\end{align*}
such that the implicit equations defining the symplectic Lorentz map are satisfied if and only if
\begin{align*}
Z_1(x,v,\Delta,\overline{v},\hbar) = Z_2(x,v,\Delta,\overline{v},\hbar) = 0.
\end{align*}
When $\hbar=0$, the unique solution of these equations for $\Delta$ and $\overline{v}$ is
\begin{align*}
\Delta_0 &= (1-\exp(-\theta_0\,\mathbb{J}) )\mathbb{J}\,(v-X_H(x)),\\
\overline{v}_0 & = X_H(x) - \exp(-\theta_0\,\mathbb{J}) [v - X_H(x)].
\end{align*}
Moreover, for each $(x,v)$, the $(\Delta,\overline{v})$ derivative of $(Z_1,Z_2)^T$ at $(x,v,\Delta_0,\overline{v}_0,0)$ is given by
\begin{align*}
\begin{pmatrix}
D_{\Delta}Z_1(x,v,\Delta_0,\overline{v}_0) & D_{\overline{v}}Z_1(x,v,\Delta_0,\overline{v}_0)\\
D_{\Delta}Z_2(x,v,\Delta_0,\overline{v}_0) & D_{\overline{v}}Z_2(x,v,\Delta_0,\overline{v}_0)
\end{pmatrix} = \begin{pmatrix}
-\mathbb{J}(x) & 1\\
\frac{1}{2}\left(\frac{\sin\theta_0}{1-\cos\theta_0}\right) & \frac{1}{2}
\end{pmatrix},
\end{align*}
which is invertible. The implicit function theorem therefore implies there is a unique pair of smooth functions $\Delta(x,v,\hbar),\overline{v}(x,v,\hbar)$ defined in an open neighborhood of $\{(x,v,\hbar)\mid \hbar=0\}\subset\mathbb{R}^n\times\mathbb{R}^n\times\mathbb{R}$ that satisfy the equations \[Z_1(x,v,\Delta(x,v,\hbar),\overline{v}(x,v,\hbar))=Z_2(x,v,\Delta(x,v,\hbar),\overline{v}(x,v,\hbar))=0.\] Since $\Delta$ is related to $\overline{q}$ by
\[
\overline{q} = x +\frac{1}{2} \hbar\,X_H(x) + \frac{1}{2}\hbar^2\,\Delta(x,v,\hbar),
\] 
another simple application of the implicit function theorem establishes existence and smoothness of the symplectic Lorentz map $\Psi_{\hbar}:(q,v)\mapsto (\overline{q},\overline{v})$. We have also shown that $\Psi_0$ has the desired form $(q,v)\mapsto (q,\overline{v}_0)$. Symplecticity of $\Psi_{\hbar}$ now follows immediately from applying the exterior derivative to \eqref{proof_symp_condition}.
\end{proof}

\begin{corollary}
The symplectic Lorentz map is a presymplectic nearly-periodic map. It is non-resonant provided $\theta_0/2\pi\nin \mathbb{Q}$.
\end{corollary}

So much for constructing a nearly-periodic map with the desired roto-rate and integral invariant. Now we must establish the precise sense in which the symplectic Lorentz map $\Psi_{\hbar}$, which is \emph{a priori} a mapping $TQ\rightarrow TQ$, functions as a consistent numerical integrator for the Hamiltonian system $X_H$ on $Q$. The first hint as to how this might work is that the limit map $\Psi_0$ admits a manifold of fixed points given by $\Gamma_0 = \{(x,v)\in TQ\mid v = X_H(q)\}$. This limiting invariant manifold, being the graph of $X_H$, is manifestly diffeomorphic to $Q$. Thus, if $\Gamma_0$ can be continued to an invariant manifold $\Gamma_{\hbar}$ for $\Psi_{\hbar}$ with $\hbar\neq 0$, we would automatically obtain dynamics on $Q$ than could be compared with those of $X_H$ by restricting $\Psi_{\hbar}$ to $\Gamma_{\hbar}$. 

Unfortunately, $\Gamma_0$ is unlikely to continue as a true invariant manifold since each fixed point on $\Gamma_0$ is of elliptic type. Instead, we can obtain the following weaker result. Roughly speaking, it says that there is a unique invariant continuation of $\Gamma_0$ at the level of formal power series in $\hbar$.

\begin{proposition}
Denote the components of the symplectic Lorentz map as $\Psi_{\hbar} = (\Psi_{\hbar}^q,\Psi_{\hbar}^v):(q,v)\mapsto (\overline{q},\overline{v})$. Assume $\theta_0\neq 0\text{ mod }2\pi$. There exists a formal power series 
\begin{align*}
v^*_{\hbar}(q) = X_H(q) + \hbar\,v^*_1(q) + \hbar^2\,v^*_2(q) + \dots
\end{align*}
with vector field coefficients such that
\begin{align}
\Psi_{\hbar}^v(q,v^*_{\hbar}(q)) = v^*_{\hbar}(\psi_{\hbar}(q)),\label{discrete_invariance}
\end{align}
where
\[
\psi_{\hbar}(q) = \Psi_{\hbar}^q(q,v^*_{\hbar}(q)).
\]
\end{proposition}

\begin{proof}
Expanding the condition \eqref{discrete_invariance} in powers of $\hbar$ leads to an infinite sequence of constraints that the formal power series $v^*_{\hbar}$ must obey. Simultaneous satisfaction of each constraint in the sequence is equivalent to \eqref{discrete_invariance}. The first two constraints are given explicitly by
\begin{align*}
\Psi_0^{v}(q,X_H(q)) &= X_H(q),\\
\Psi_1^v(q,X_H(q)) + D_v\Psi_0^v(q,X_H(q))[v^*_1(q)] &= v^*_1(q) + DX_H(q)[\Psi_{1}^{q}(q,X_H(q))],
\end{align*}
where we have introduced the formal series expansions 
\begin{align*}
\Psi_{\hbar}^q = \Psi^q_0 + \hbar\,\Psi^q_1 + \dots,\quad \Psi_{\hbar}^v = \Psi_{0}^v + \hbar\,\Psi_{1}^v+\dots
\end{align*}
and used $\Psi^q_0(q,v)=q$. Glancing at \eqref{limit_symplectic_lorentz} reveals that the first of these equations is automatically satisfied. The second equation can be interpreted as an algebraic equation constraining the form of $v^*_1$. In fact, since the linear map 
\begin{gather*}
L(q) = D_v\Psi_0^v(q,X_H(q)) - \text{id},\\
\delta v\mapsto \exp(-\theta_0\,\mathbb{J}(q))[\delta v] - \delta v,
\end{gather*} 
is invertible whenever $\theta_0\neq 0\text{ mod }2\pi$, $v^*_1$ is determined uniquely by the formula 
\begin{align*}
v^*_1(q) = L(q)^{-1}[DX_H(q)[\Psi_{1}^{q}(q,X_H(q))] - \Psi_1^v(q,X_H(q)) ].
\end{align*}
More generally, the $n^{\text{th}}$ equation in the sequence has the form
\begin{align*}
L(q)[v_n^*(q)] = s_n(q),
\end{align*}
where $s_n(q)$ depends only on coefficients of the power series expansion for $\Psi_{\hbar}$ and coefficients $v^*_k$ with $k<n$. Invertibility of $L(q)$ therefore implies that there is a unique formula for $v^*_n$ for each $n$. The formal power series $v^*_{\hbar}$ defined in this manner satisfies \eqref{discrete_invariance} by construction.
\end{proof}

So while we do not obtain a genuine invariant manifold diffeomorphic to $Q$, we do obtain a family of approximate invariant manifolds diffeomorphic to $Q$ given by truncations of the formal power series $v^*_{\hbar}$. Using arguments comparable to those presented in \cite{Burby_Hirvijoki_2021_JMP}, it is possible to show that truncations of $v^*_{\hbar}$ may be constructed so their graphs agree with the zero level set of the adiabatic invariant $\mu_{\hbar}$ for $\Psi_{\hbar}$ to any desired order in $\hbar$. Adiabatic invariance of $\mu_{\hbar}$ can then be used to prove the existence of manifolds $\Gamma_{\hbar}^{(n)}$ close to $\Gamma_0$ with the following schematic normal stability property:
\begin{itemize}
\item For each $N>0$ and $(q,v)$ within $\hbar^{\alpha(n)}$ of $\Gamma_{\hbar}^{(n)}$ the point $\Psi_{\hbar}^{k}(q,v)$ remains within $\hbar^{\beta(n)}$ of $\Gamma_{\hbar}^{(n)}$ for $k=O(\epsilon^{-N})$. Here $\alpha$ and $\beta$ are monotone increasing functions of $n$.
\end{itemize}
We will not attempt to prove such a result in full generality here. However, since we have already determined the form of the leading term in the adiabatic invariant series (in this case $\mu_{\hbar} = \hbar^4\,\mu_4 + O(\hbar^5)$), we can prove a special case of the general result without much effort. First we we establish the timescale over which $\mu_4$ is well-conserved.

\begin{definition}
Given a compact set $C\subset TQ$, a point $(q,v)$ is \textbf{positively-contained} if $\Psi_{\hbar}^k(q,v)\in C$ for all non-negative integers $k$.
\end{definition}

\begin{proposition}\label{adiabatic_invariance_zeroth_order}
For each $N>0$ and compact set $C\subset TQ$ there is a positive, $\hbar$-independent constant $\mathcal{M}$ such that
\begin{align*}
|\mu_4(\Psi_{\hbar}^k(q,v))-\mu_4(q,v)|\leq \mathcal{M}\,\hbar,\quad k\in [0,k^*(\hbar,N)],
\end{align*}
whenever $(q,v)$ is positively-contained in $C$. Here $k^*(\hbar,N) = O(\hbar^{-N})$. 
\end{proposition}
\begin{proof}
First we will obtain a useful estimate for the degree of conservation of an arbitrary truncation of the adiabatic invariant series. Let $ \overline{\mu}_{\hbar} = \hbar^{-4}\mu_{\hbar} = \mu_4 + \hbar\,\mu_5 + \dots$ denote the reduced adiabatic invariant for the symplectic Lorentz map. Define $\overline{\mu}_i = \mu_{i+4}$ and let $\overline{\mu}_{\hbar}^{(N)} = \sum_{i=0}^{N}\overline{\mu}_i\,\hbar^{i}$. Since $\overline{\mu}_{\hbar}^{(N)} = \overline{\mu}_{\hbar} + O(\hbar^{N+1})$ and $\overline{\mu}_{\hbar}$ is $\Psi_{\hbar}$-invariant to all orders in $\hbar$, there is a constant $\mathcal{M}_N$, depending on both $C$ and $N$, such that
\begin{align*}
\forall (q,v)\in C,\,|\overline{\mu}_{\hbar}^{(N)}(\Psi_{\hbar}(q,v)) - \overline{\mu}_{\hbar}^{(N)}(q,v) | \leq \mathcal{M}_N\,\hbar^{N+1} .
\end{align*}
For positively-contained $(q,v)$ we may apply this formula repeatedly to obtain an estimate for the change in $\overline{\mu}_{\hbar}^{(N)}$ after $k$ positive timesteps,
\begin{align}
|\overline{\mu}_{\hbar}^{(N)}(\Psi_{\hbar}^k(q,v)) - \overline{\mu}_{\hbar}^{(N)}(q,v) |&\leq \mathcal{M}_{N}\,\hbar^{N+1} + |\overline{\mu}_{\hbar}^{(N)}(\Psi_{\hbar}^{k-1}(q,v)) - \overline{\mu}_{\hbar}^{(N)}(q,v) |\nonumber\\
&\leq (1+k)\mathcal{M}_N\,\hbar^{N+1}.\label{k_step_bound_N}
\end{align}

Next, we draw implications from the previous inequality together with a bound on the difference between $\overline{\mu}_{\hbar}^{(0)} = \mu_4$ and $\overline{\mu}_{\hbar}^{(N)}$. There must be another positive constant $\mathcal{M}^\prime_{N}$, depending on both $C$ and $N$, such that
\begin{align*}
\forall (q,v)\in C,\,|\overline{\mu}_{\hbar}^{(0)}(q,v) - \overline{\mu}_{\hbar}^{(N)}(q,v) | \leq \mathcal{M}^\prime_N\,\hbar.
\end{align*}
In light of the inequality \eqref{k_step_bound_N}, this implies that for each positively-contained $(q,v)$ the change in $\overline{\mu}_{\hbar}^{(0)}$ after $k$ positive timesteps is at most
\begin{align*}
|\overline{\mu}_{\hbar}^{(0)}(\Psi_{\hbar}^k(q,v)) - \overline{\mu}_{\hbar}^{(0)}(q,v) |&\leq |\overline{\mu}_{\hbar}^{(0)}(\Psi_{\hbar}^k(q,v)) - \overline{\mu}_{\hbar}^{(N)}(\Psi^k_{\hbar}(q,v)) | + | \overline{\mu}_{\hbar}^{(N)}(\Psi^k_{\hbar}(q,v))- \overline{\mu}_{\hbar}^{(0)}(q,v) |\\
&\leq \mathcal{M}^\prime_{N}\,\hbar + | \overline{\mu}_{\hbar}^{(N)}(\Psi_{\hbar}^k(q,v))- \overline{\mu}_{\hbar}^{(0)}(q,v) |\\
&\leq \mathcal{M}^\prime_{N}\,\hbar +  | \overline{\mu}_{\hbar}^{(N)}(\Psi^k_{\hbar}(q,v))- \overline{\mu}_{\hbar}^{(N)}(q,v) | +  |\overline{\mu}_{\hbar}^{(N)}(q,v)- \overline{\mu}_{\hbar}^{(0)}(q,v) |\\
&\leq 2\mathcal{M}^\prime_{N}\,\hbar + (1+k)\mathcal{M}_N\,\hbar^{N+1}.
\end{align*}
Apparently, the change in $\mu_4 = \overline{\mu}_{\hbar}^{(0)}$ is at most $O(\hbar)$ as long as $k\,\hbar^{N+1}=O(\hbar)$. We therefore obtain the desired inequality with $k^*(\hbar,N) = \lfloor{\hbar}^{-N}\rfloor -1$.

\end{proof}

Using this result together with the explicit form of $\mu_4$, we now easily obtain the following normal stability result for the almost invariant set given by the graph of $X_H$.

\begin{proposition}\label{normal_stability_estimate_zeroth}
Let $C\subset TQ$ be a compact set and set $(q^k,v^k)=\Psi_{\hbar}^k(q,v)$ for any $(q,v)\in TQ$. Let $|\cdot|$ denote the velocity norm provided by the metric tensor $g$. For each $N > 0$, $V_0>0$, and positively-contained $(q,v)\in C$ that satisfies
\begin{align*}
|v-X_H(q)|_q<V_0\,\sqrt{\hbar},
\end{align*}
there is a positive constant $V_1$ such that
\begin{align*}
|v^k-X_H(q^k)|_{q^k}\leq V_1\,\sqrt{\hbar}
\end{align*}
for all $k\in[0,k^*(\hbar,N)]$. Here, $k^*(\hbar,N) = O(\hbar^{-N})$.
\end{proposition}

\begin{proof}
Let $(q,v)\in C$ be positively-contained and suppose  $|v-X_H(q)|_q<V_0\,\hbar$. By Proposition \ref{adiabatic_invariance_zeroth_order}, we have
\begin{align*}
|\mu_4(\Psi_{\hbar}^k(q,v))-\mu_4(q,v)|\leq \mathcal{M}\,\hbar,
\end{align*}
for some $N$-dependent constant $\mathcal{M}$ and $k\in[0,k^*(\hbar,N)]$. But since $\mu_4(q,v) =\frac{1}{2} |v-X_H(q)|_q^2$ we can apply this inequality to obtain
\begin{align*}
|v^k-X_H(q^k)|_{q^k}^2 &\leq ||v^k-X_H(q^k)|_{q^k}^2 - |v-X_H(q)|_q^2| + |v-X_H(q)|_q^2\\
&\leq 2|\mu_4(\Psi_{\hbar}^k(q,v))-\mu_4(q,v)| + V_0^2\,\hbar\\
&\leq 2\mathcal{M}\,\hbar + V_0^2\,\hbar.
\end{align*}
%\begin{align*}
%\left|g_{q^k}(v^k-X_H(q^k),v^k-X_H(q^k)) - g_q(v-X_H(q),v-X_H(q))\right|\leq 2 \mathcal{M}\,\hbar
%\end{align*}
Taking a square root gives the desired result.
\end{proof}

In the above sense, the graph of $X_H$ behaves much like a true invariant set over very large time intervals. Of course, the invariance need not be exact, but may include oscillations around the graph of amplitude $\sqrt{\hbar}$. The amplitude of these oscillations can be reduced by considering manifolds that better approximate the zero level set of $\mu_{\hbar}$, but, as mentioned earlier, we will not pursue this matter further in this article.

To complete the picture of how the symplectic Lorentz map may be used as an integrator $X_H$ on $Q$, we will now describe the precise sense in which the map's dynamics approximate the $H$-flow. We start with a simple estimate that says the $q$-component of the symplectic Lorentz map approximates the time-$\hbar$ flow of $X_H$ on $Q$ with an $O(\hbar^{5/2})$ error, \emph{provided} the map is applied in an $O(\hbar^{1/2})$ neighborhood of the graph $\{v=X_H(q)\}$.

\begin{proposition}\label{local_accuracy_slm}
Let $(q,v_{\hbar})$ be a smooth $\hbar$-dependent point in $TQ$ with $v_{\hbar} = X_H(q) + O(\hbar^{1/2})$. The mapped point $(\overline{q},\overline{v})=\Psi_{\hbar}(q,v_{\hbar})$ satisfies
\begin{align*}
\overline{q}&= q + \hbar\,X_H(q) +\frac{1}{2}\hbar^2\,DX_H(q)[X_H(q)] + O(\hbar^{5/2}),\\
\overline{v}&=X_H(\overline{q})+O(\hbar^{1/2}).
\end{align*}
\end{proposition}

\begin{proof}
In the proof of Proposition \ref{slm_existence}, we already established
\begin{align*}
\overline{q} &= x +\frac{1}{2} \hbar\,X_H(x) + \frac{1}{2}\hbar^2\,\Delta(x,v_{\hbar},\hbar),\\
\overline{v} & = X_H(x) - \exp(-\theta_0\,\mathbb{J}_x) [v_{\hbar} - X_H(x)]+O(\hbar),
\end{align*}
and 
\begin{align*}
\Delta(x,v_{\hbar},0) &= (1-\exp(-\theta_0\,\mathbb{J}) )\mathbb{J}\,(v_{\hbar}-X_H(x)),
\end{align*}
where $x= \overline{q}/2 + q/2$. Implicit differentiation of these formulas together with Taylor's theorem with remainder therefore implies
\begin{align*}
\overline{q} & = q + \hbar\,X_H(q)+\frac{1}{2}\hbar^2\,DX_H(q)[X_H(q)] + \hbar^2\,(1-\exp(-\theta_0\,\mathbb{J}_q) )\mathbb{J}_q\,(v_{\hbar}-X_H(q)) + O(\hbar^3),\\
\overline{v} & = X_H(q) - \exp(-\theta_0\,\mathbb{J}_q) [v_{\hbar} - X_H(q)]+O(\hbar)
\end{align*}
The desired result now follows immediately from $v_{\hbar}-X_H(q) = O(\hbar^{1/2})$.
\end{proof}

Combining this result with our earlier estimate of the normal stability timescale for $\{v=X_H(q)\}$ in Proposition \ref{normal_stability_estimate_zeroth} finally allows us to conclude that the $q$-component of the symplectic Lorentz map provides a persistent approximation of the $H$-flow over very large time intervals provided initial conditions are chosen close enough to the graph $\{v=X_H(q)\}$.

\begin{corollary}[Persistent approximation property]
Let $C$ be a compact set and let $(q,v_{\hbar})\in C$ be a smooth $\hbar$-dependent point in $C$ that is positively-contained for each $\hbar$. Also assume $v_{\hbar} = X_H(q) + O(\hbar^{1/2})$.  For each $N>0$ there is an integer $k^*(\hbar,N)=O(\hbar^{-N})$ such that 
\begin{align*}
q^{k+1}&= q^k + \hbar\,X_H(q^k) +\frac{1}{2}\hbar^2\,DX_H(q^k)[X_H(q^k)] + O(\hbar^{5/2}),\\
v^{k+1}&=X_H(q^{k+1})+O(\hbar^{1/2}),
\end{align*}
for each $k\in[0,k^*(\hbar,N)]$. Here, $(q^k,v^k) = \Psi_{\hbar}^k(q,v_{\hbar})$.
\end{corollary}

\begin{proof}
Proposition \ref{normal_stability_estimate_zeroth} ensures that the iterates $(q^k,v^k)$ remain within $O(\hbar^{1/2})$ of $\{v=X_H(q)\}$ for $k$ in the desired range. Thus, Proposition \ref{local_accuracy_slm} applies to each iterate individually, which is precisely the desired result.
\end{proof}

In summary, we have established the following remarkable properties of the symplectic Lorentz map $\Psi_{\hbar}$.
\begin{enumerate}
\item It is symplectic on $TQ$, when $TQ$ is endowed with the magnetic symplectic form $\Omega^*_{\hbar} = \pi^*\omega + \hbar^2\,\Omega$.
\item Its $q$-component provides an approximation of the time-$\hbar$ flow of $X_H$ with $O(\hbar^{5/2})$ local truncation error when applied to points in $TQ$ within $O(\hbar^{1/2})$ of $\{v=X_H(q)\}$.
\item If an initial condition is chosen to lie within $\hbar^{1/2}$ of $\{v=X_H(q)\}$ then it will remain within $\hbar^{1/2}$ of the same set for a number of iterations that scales like $\hbar^{-N}$ for any $N$.
\end{enumerate}

\section{Examples\label{examples_section}}
\subsection{Hidden-variable Newtonian gravity\label{hidden_gravity_section}}

In this section, we will use nearly-periodic maps to construct a discrete-time model of Newtonian gravitation where the gravitational constant has a dynamical origin. Let $M = \mathbb{R}\times\mathbb{R}\times\mathbb{R}^{d}\times\mathbb{R}^d\ni (q,p,\bm{Q},\bm{P})$ and set $\Omega_\gamma = dq\wedge dp + \sum_{i=1}^d dQ^i\wedge dP_i$. Let $V,W:\mathbb{R}^d\rightarrow \mathbb{R}$ be smooth functions. The Hamiltonian 
\begin{align}
H_\epsilon(q,p,\bm{Q},\bm{P}) & = \frac{1}{2}\,(p^2 + q^2) + \epsilon\,\bigg(\frac{1}{2}|\bm{P}|^2+V(\bm{Q})  +q^2\,W(\bm{Q})\bigg)
\end{align}
defines a continuous-time nearly-periodic Hamiltonian system with equations of motion
\begin{align}
\dot{p}& = -q - \epsilon\,2\,q\,W(\bm{Q})\label{Bex_fs_1}\\
\dot{q}& = p\\
\dot{\bm{P}}& =  -{\epsilon}\,\partial_{\bm{Q}}V- \epsilon\, q^2\, \partial_{\bm{Q}}W\\
\dot{\bm{Q}}& = \epsilon\,\bm{P}.\label{Bex_fs_4}
\end{align}
 The angular frequency function is $\omega_0 = 1$, the limiting roto-rate is $R_0 = -q\,\partial_p + p\,\partial_q$, and the corresponding $U(1)$-action is $\Phi_\theta(q,p,\bm{Q},\bm{P}) = (\cos\theta\,q + \sin\theta\,p,\cos\theta\,p-\sin\theta\,q,\bm{Q},\bm{P})$. When $\epsilon=0$, the system's flow is $F_t(q,p,\bm{Q},\bm{P}) = \Phi_t(q,p,\bm{Q},\bm{P})$. Intuitively, the $(q,p)$ variables correspond to a fast oscillator that couples nonlinearly to a mechanical system parameterized by $(\bm{Q},\bm{P})$. The averaged Hamiltonian for the coupled system is 
 \begin{align}
 \frac{1}{2\pi}\int_0^{2\pi}\Phi_\theta^*H_\epsilon\,d\theta = \mu_0+\epsilon\,\bigg(\frac{1}{2}|\bm{P}|^2+V(\bm{Q})  +\mu_0\,W(\bm{Q})\bigg),
 \end{align}
 where $\mu_0 = \frac{1}{2}(p^2+q^2)$ is the leading-order adiabatic invariant. We therefore expect the slow variables $(\bm{Q},\bm{P})$ to behave like a particle in $d$-dimensional space subject to the effective potential $V(\bm{Q})+ \mu_0\,W(\bm{Q})$.
 
 We will construct a Hamiltonian non-resonant nearly-periodic map that accurately simulates the slow dynamics for this system while ``stepping over" the shortest scale $2\pi/\omega_0\sim 1$. If $h\in\mathbb{R}$ denotes the temporal step size, these requirements translate into symbols as $1\ll h\ll \epsilon^{-1}$. Upon introducing the parameters $\delta = 1/h$, $\hbar = \epsilon\,h$, and $\gamma = (\hbar,\delta)$, we may state our requirement equivalently as $|\gamma|\ll 1$. Our construction will now proceed using the method of mixed-variable generating functions.
 
 The exact Type I generating function for this problem can be characterized by Jacobi's solution of the Hamilton--Jacobi equation, which is given by
 \begin{align}
     S(q,\bm{Q},\overline{q},\overline{\bm{Q}})=\int_0^h \left[ p(t)\dot{q}(t) + \bm{P}(t)\dot{\bm{Q}}(t)-H(q(t),\bm{Q}(t),p(t),\bm{P}(t))\right] dt,
 \end{align}
 where $(q(t),\bm{Q}(t),p(t),\bm{P}(t))$ satisfies Hamilton's equations, and the boundary conditions $q(0)=q$, $\bm{Q}(0)=\bm{Q}$, $q(h)=\overline{q}$, $\bm{Q}(h)=\overline{\bm{Q}}$. In the setting of variational integrators~\cite{MaWe2001}, this is referred to as the exact discrete Lagrangian, and there are also exact discrete Hamiltonians~\cite{LeZh2009} corresponding to Type II and Type III generating functions. One possible way to construct a computable approximation of the exact Type I generating function is to observe that it can also be expressed as,
 \begin{align}
     S(q,\bm{Q},\overline{q},\overline{\bm{Q}})=\ext_{\substack{(q, p, \bm{Q},\bm{P}) \in
C^2([0, h],M)\\q(0)=q, q(h)=\overline{q},\\
\bm{Q}(0)=\bm{Q}, \bm{Q}(h)=\overline{\bm{Q}}}}\int_0^h \left[ p(t)\dot{q}(t) + \bm{P}(t)\dot{\bm{Q}}(t)-H(q(t),\bm{Q}(t),p(t),\bm{P}(t))\right] dt,
 \end{align}
Then, one can construct a computable approximation by replacing the infinite-dimensional function space $C^2([0, h],M)$ with a finite-dimensional subspace, and replacing the integral with a numerical quadrature formula, which yields a Galerkin discrete Lagrangian. Under a number of technical assumptions, the resulting variational integrators $\Gamma$-converge to the exact flow map~\cite{MuOr2004}, and a quasi-optimality result~\cite{HaLe2015} implies that the rate of convergence is related to the best approximation properties of the finite-dimensional function space used to construct the Galerkin discrete Lagrangian. In general, this means that a good integrator can be constructed by choosing a finite-dimensional function space that is rich enough to approximate the exact solutions well, and using a quadrature rule that is accurate for that choice of function space.

This might entail augmenting the function space with the solution of the fast dynamics when the slow variables are frozen, and then using a quadrature rule that is well-adapted to highly oscillatory integrals, like Filon quadrature~\cite{IsNo2004}. In this case, the problem exhibits a fast-slow structure that lends itself to a hybrid approximation. We exploit the time-scale separation to approximate the fast variables of the dynamics $(q(t),p(t))$ by the exact solution of the $\epsilon=0$ limiting system, 
\begin{align*}
\dot{p} = -q, \quad
\dot{q} = p, \quad
\dot{\bm{P}} =0, \quad
\dot{\bm{Q}} = 0.
\end{align*}
where the slow variables $(\bm{Q}(t),\bm{P}(t))$ are frozen, which leads to a sinusoidal solution for $(q,p)$. Furthermore, because the timestep $h$ is assumed to be large enough that the fast variables perform many revolutions in that time, we anti-alias the fast dynamics by replacing the revolutions by just the fractional part of the revolutions, which we denote by $\theta_0$, and which is assumed to be some irrational multiple of $2\pi$, so that the invariant distribution remains the same. The component of the action integral associated with the fast variables can be evaluated analytically in this case. As for the slow variables, we adopt an approach that can be used to derive the implicit midpoint rule, which is a symplectic integrator for Hamiltonian systems. This involves approximating the solution space by linear functions, so $\bm{Q}(t)$ is uniquely determined by the boundary conditions, and approximating the integral by the midpoint rule. The use of mixed quadrature approximations of the action integral was the basis for implicit-explicit variational integrators for fast-slow systems~\cite{StGr2009}. 

Note that $\omega_0=1$, then for the $(q,p)$ dynamics to have a $\theta_0$ rotation in time $h$, the solution is given by,
\[
    \begin{bmatrix}
      q(t) \\
      p(t)
    \end{bmatrix} =
    \begin{bmatrix}
      \cos (t) & \sin (t)\\
      -\sin (t) & \cos (t)
    \end{bmatrix}
    \begin{bmatrix}
      q \\ p
    \end{bmatrix},
\]
where $\theta_0/h=1$. Since $q(t)=\sin (t) p + \cos (t) q$, then $\overline{q}=q(h)=\sin (\theta_0) p + \cos (\theta_0) q$, and hence, $p=\frac{\overline{q}-\cos (\theta_0) q}{\sin (\theta_0)}$. Therefore, the $(q,p)$ dynamics, expressed in terms of the boundary data is given by,
\[
\begin{bmatrix}
      q(t) \\
      p(t)
    \end{bmatrix} =
    \begin{bmatrix}
      \cos (t) & \sin (t)\\
      -\sin (t) & \cos (t)
    \end{bmatrix}
    \begin{bmatrix}
    q\\
      \frac{\overline{q}-\cos (\theta_0) q}{\sin (\theta_0)}
    \end{bmatrix}
    =\begin{bmatrix}
      \cos (t) q+\sin(t) \frac{\overline{q}-\cos (\theta_0) q}{\sin (\theta_0)}\\
    -\sin(t) q+\cos(t) \frac{\overline{q}-\cos (\theta_0) q}{\sin (\theta_0)} 
    \end{bmatrix}.
\]
For the slow degrees of freedom, we consider a linear interpolant in time, $\bm{Q}(t)=\bm{Q}+\frac{\overline{\bm{Q}}-\bm{Q}}{h}t$, from which the momentum becomes  $\bm{P}(t)=\frac{1}{\epsilon}\dot{\bm{Q}}(t)=\frac{1}{\epsilon}\frac{\overline{\bm{Q}}-\bm{Q}}{h}$. Collecting all of these, we have
%\begin{align*}
%    q(t) &= \sin(t) \frac{\overline{q}-\cos (\theta_0) q}{\sin (\theta_0)} + \cos (t) q\\
%    \dot{q}(t) & = \cos(t) \frac{\overline{q}-\cos (\theta_0) q}{\sin (\theta_0)} -\sin (t) q\\
%    p(t)&=\cos(t) \frac{\overline{q}-\cos (\theta_0) q}{\sin (\theta_0)} -\sin(t) q\\
%    \bm{Q}(t)&=\bm{Q}+\frac{\overline{\bm{Q}}-\bm{Q}}{h}t\\
%    \dot{\bm{Q}}(t)&= \frac{\overline{\bm{Q}}-\bm{Q}}{h}\\
%    \bm{P}(t)&=\frac{1}{\epsilon}\frac{\overline{\bm{Q}}-\bm{Q}}{h}.
%\end{align*}
\begin{alignat*}{5}
    q(t) &= \cos (t) q+\sin(t) \frac{\overline{q}-\cos (\theta_0) q}{\sin (\theta_0)},
    &\qquad\qquad
    \bm{Q}(t)&=\bm{Q}+\frac{\overline{\bm{Q}}-\bm{Q}}{h}t,\\
    \dot{q}(t) & = -\sin (t) q+\cos(t) \frac{\overline{q}-\cos (\theta_0) q}{\sin (\theta_0)},
    &\qquad\qquad
    \dot{\bm{Q}}(t)&= \frac{\overline{\bm{Q}}-\bm{Q}}{h},\\
    p(t)&=-\sin(t) q+\cos(t) \frac{\overline{q}-\cos (\theta_0) q}{\sin (\theta_0)} ,
    &\qquad\qquad
    \bm{P}(t)&=\frac{1}{\epsilon}\frac{\overline{\bm{Q}}-\bm{Q}}{h}.
\end{alignat*}
With these, we are now ready to approximate the exact Type I generating function,
\begin{align*}
 S(q,\bm{Q},\overline{q},\overline{\bm{Q}})&=\int_0^h \left[ p(t)\dot{q}(t) + \bm{P}(t)\dot{\bm{Q}}(t)-H(q(t),\bm{Q}(t),p(t),\bm{P}(t))\right] dt\\
 &=\int_0^h \Bigg[ p(t)\dot{q}(t) -\frac{1}{2}\,(p(t)^2 + q(t)^2)\\
 &\qquad\qquad + \bm{P}(t)\dot{\bm{Q}}(t)-\epsilon\,\bigg(\frac{1}{2}|\bm{P}(t)|^2+V(\bm{Q}(t))  +q(t)^2\,W(\bm{Q}(t))\bigg)\Bigg] dt,
\end{align*}
where we observe that the integral involving the $(q,p)$ terms have the form of an trigonometric integral, which can be evaluated analytically,
\begin{align*}
    &\int_0^h \left[p(t)\dot{q}(t)-\frac{1}{2}\,(p^2 + q^2)\right]\\
    &\qquad= \int_0^h \left[(p^2\cos^2(t)-2pq\cos(t)\sin(t)+q^2\sin^2(t))-\frac{1}{2}(p^2+q^2)\right]dt\\
    &\qquad=\int_0^h \left[\left(p^2\left(\frac{1+\cos (2t)}{2}\right)-pq\sin(2t)+q^2\left(\frac{1-\cos (2t)}{2}\right)\right)-\frac{1}{2}(p^2+q^2)\right]dt\\
    &\qquad=\int_0^h \left[\left(\frac{1}{2}(p^2-q^2)\cos (2t)-pq\sin(2t)\right)\right]dt\\
    &\qquad=\left.\frac{1}{4}(p^2-q^2)\sin(2t)+\frac{1}{2}pq\cos(2t)\right|_0^{\theta_0}\\
    &\qquad=\frac{\cos(\theta_0)\tfrac{1}{2}\overline{q}^2+\cos(\theta_0)\tfrac{1}{2}q^2-q\overline{q}}{\sin(\theta_0)}
\end{align*}
where we made use of the trigometric double angle formulas, $p=\frac{\overline{q}-\cos (\theta_0) q}{\sin (\theta_0)}$, and $h=\theta_0$. It is easy to verify that this generating function generates a $\theta_0$ rotation in the $(q,p)$ variables.

We approximate the integral involving the $(\bm{Q},\bm{P})$ terms with the midpoint rule,
 \begin{align*}
     &\int_0^h \left[\bm{P}(t)\dot{\bm{Q}}(t) - \epsilon\,\bigg(\frac{1}{2}|\bm{P}(t)|^2+V(\bm{Q}(t))  +q(t)^2\,W(\bm{Q}(t))\bigg)\right]dt\\
     &\qquad\approx
     h\left[\frac{1}{\epsilon}\left(\frac{\overline{\bm{Q}}-\bm{Q}}{h}\right)^2-\epsilon\left(\frac{1}{2}\frac{1}{\epsilon^2}\left(\frac{\overline{\bm{Q}}-\bm{Q}}{h}\right)^2+V\left(\frac{\bm{Q}+\overline{\bm{Q}}}{2}\right)+q\left(\frac{h}{2}\right)^2W\left(\frac{\bm{Q}+\overline{\bm{Q}}}{2}\right)\right)\right]\\
     &\qquad=  h\left[\frac{1}{2\epsilon}\left(\frac{\overline{\bm{Q}}-\bm{Q}}{h}\right)^2-\epsilon\left(V\left(\frac{\bm{Q}+\overline{\bm{Q}}}{2}\right)+q\left(\frac{h}{2}\right)^2W\left(\frac{\bm{Q}+\overline{\bm{Q}}}{2}\right)\right)\right].
 \end{align*}
 By using $\hbar=\epsilon h$, replacing $q\left(\frac{h}{2}\right)$ with $q(0)=q$, and combining it with the first term coming from the fast dynamics, we obtain the following Type I generating function,
 \begin{align}
     S_\gamma(q,\bm{Q},\overline{q},\overline{\bm{Q}}) & = \frac{\cos\theta_0\,\tfrac{1}{2}q^2 + \cos\theta_0\tfrac{1}{2}\overline{q}^2 - q\overline{q}}{\sin\theta_0} +\frac{1}{2}\frac{|\overline{\bm{Q}}-{\bm{Q}}|^2}{\hbar}\nonumber\\
     &\qquad-\hbar\,V\left(\frac{\bm{Q}+\overline{\bm{Q}}}{2}\right)- \hbar\, q^2 W\left(\frac{\bm{Q}+\overline{\bm{Q}}}{2}\right),
 \end{align}
 where $\theta_0$ is some irrational multiple of $2\pi$. The implicit relations
 \begin{gather*}
    \overline{p}=\partial_{\overline{q}}S_\gamma,\quad p = -\partial_{q}S_\gamma,\quad\overline{\bm{P}} = \partial_{\overline{\bm{Q}}}S_\gamma,\quad \bm{P} = -\partial_{\bm{Q}}S_\gamma,\label{Bgen_defs}
 \end{gather*}
 define a $\gamma$-dependent symplectic map $F_\gamma:(q,p,\bm{Q},\bm{P})\mapsto (\overline{q},\overline{p},\overline{\bm{Q}},\overline{\bm{P}})$. For small $\gamma$, we claim this map accurately captures the averaged dynamics of the slow variables in the system \eqref{Bex_fs_1}-\eqref{Bex_fs_4} and preserves the adiabatic invariant $\mu_0$ over very large time intervals. To show this, we first compute the derivatives in \eqref{Bgen_defs} to explicitly write the defining equations for $F_\gamma$ as
 \begin{align}
     \overline{p} &= \frac{\cos\theta_0}{\sin\theta_0}\overline{q} - \frac{1}{\sin\theta_0}q, \nonumber\\
     p & = -\frac{\cos\theta_0}{\sin\theta_0}q + \frac{1}{\sin\theta_0}\overline{q} + \hbar\,2q\,W\left(\frac{\bm{Q}+\overline{\bm{Q}}}{2}\right),\nonumber\\
     \overline{\bm{P}} & = \frac{\overline{\bm{Q}}-\bm{Q}}{\hbar}-\frac{1}{2}\hbar\,V^\prime\left(\frac{\bm{Q}+\overline{\bm{Q}}}{2}\right) - \frac{1}{2}\hbar\,q^2\, W^\prime\left(\frac{\bm{Q}+\overline{\bm{Q}}}{2}\right),\nonumber\\
     \bm{P} & = \frac{\overline{\bm{Q}}-\bm{Q}}{\hbar} +\frac{1}{2}\hbar\,V^\prime\left(\frac{\bm{Q}+\overline{\bm{Q}}}{2}\right) + \frac{1}{2}\hbar \,q^2\,W^\prime\left(\frac{\bm{Q}+\overline{\bm{Q}}}{2}\right).\nonumber
 \end{align}
 The first pair of equations can be solved explicitly for $\overline{p}$ and $\overline{q}$, giving
 \begin{align}
    \label{eq:ex2_fast_slow_small_p}
     \overline{p} & = \cos\theta_0\,p - \sin\theta_0\,q - \cos\theta_0\,\hbar\,2q\,W\left(\frac{\bm{Q}+\overline{\bm{Q}}}{2}\right),\\
     \label{eq:ex2_fast_slow_small_q}
     \overline{q} & = \cos\theta_0\,q + \sin\theta_0\,p - {\sin\theta_0}\, \hbar\,2q\,W\left(\frac{\bm{Q}+\overline{\bm{Q}}}{2}\right).
\end{align}
Adding and subtracting the last $\bm{P}$ and $\overline{\bm{P}}$ equations then gives
\begin{align}
    \overline{\bm{P}} - \bm{P} & = -\hbar\, V^\prime\left(\frac{\bm{Q}+\overline{\bm{Q}}}{2}\right) - \hbar\,q^2\,W^\prime\left(\frac{\bm{Q}+\overline{\bm{Q}}}{2}\right),\label{Bim_toy1}\\
    \overline{\bm{Q}} - \bm{Q} & = \hbar \frac{\bm{P}+\overline{\bm{P}}}{2}.\label{Bim_toy2}
\end{align}
These formulas show that $F_0 = \Phi_{\theta_0}$, which implies that $F_\gamma$ comprises a non-resonant, Hamiltonian, nearly-periodic map. In particular, this map admits an all-orders adiabatic invariant $\mu = \mu_0 + O(\epsilon)$, where $\mu_0 = \tfrac{1}{2}(q^2+p^2)$. Moreover, equations \eqref{Bim_toy1}-\eqref{Bim_toy2} provide a consistent numerical scheme for the averaged dynamics of the slow variables. To see this, note that the average of $q$ in \eqref{Bim_toy1} after many iterations tends to $\mu_0$, which implies that, on average, \eqref{Bim_toy1}-\eqref{Bim_toy2} comprise the implicit midpoint scheme applied to the continuous system's averaged dynamics. Note that the relationship between the physical timestep $h$ and $\hbar$ is $h = \hbar/\epsilon$.
 
A planar $N$-body problem in Cartesian $(x,y)$-coordinates provides a convenient sandbox for testing the novel scheme \eqref{eq:ex2_fast_slow_small_p}--\eqref{Bim_toy2}. Assume two bodies, labelled by the position vectors $Q_1=(Q_{1,x},Q_{1,y})$ and $Q_2=(Q_{2,x},Q_{2,y})$ and the respective momentum vectors $P_1=(P_{1,x},P_{1,y})$ and $P_2=(P_{2,x},P_{2,y})$, to orbit an infinitely massive body at the origin. The potential $V(\bm{Q})$ is therefore
 \begin{align}
     V(Q_1,Q_2) & = -\frac{1}{|Q_1|} - \frac{1}{|Q_2|}.
\end{align}
Also assume the two bodies to interact via the additional central potential 
\begin{align}
     %\partial_{Q_1}V & = -\frac{Q_1}{|Q_1|^3},\quad \partial_{Q_2}V = -\frac{Q_2}{|Q_2|^3}\\
     W(Q_1,Q_2) & = -\frac{1}{|Q_1-Q_2|}.
\end{align}
The instantaneous value of $q^2$ therefore indicates the strength of the coupling of the two bodies via the temporal evolution of the $\epsilon$-perturbed $(q,p)$ oscillator. 

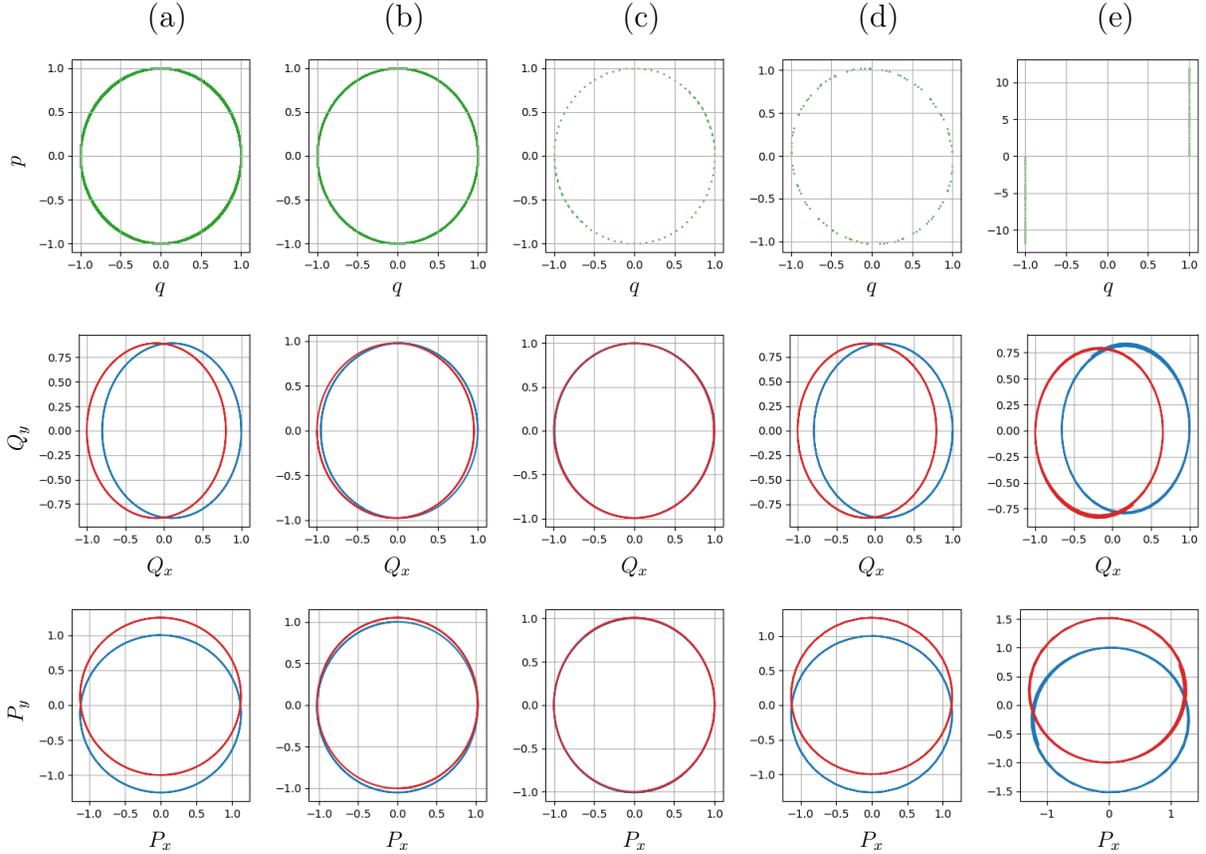
\begin{figure}[!htb]
    
    \centering
    
    \begin{tikzpicture}[image/.style = {text width=0.19\textwidth,inner sep=0pt, outer sep=0pt},node distance = 0mm and 0mm] 
        %%%%
        \node [image] (frame11)
        {\includegraphics[width=1.0\linewidth]{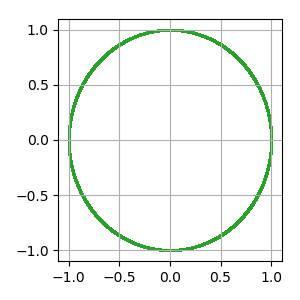}};
        
        \node [image,below=0.5cm of frame11] (frame21)
        {\includegraphics[width=1.0\linewidth]{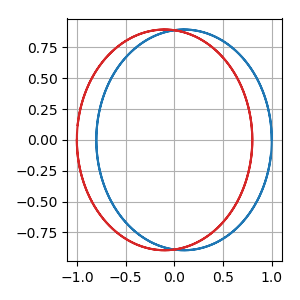}};
        
        \node [image,below=0.5cm of frame21] (frame31)
        {\includegraphics[width=1.0\linewidth]{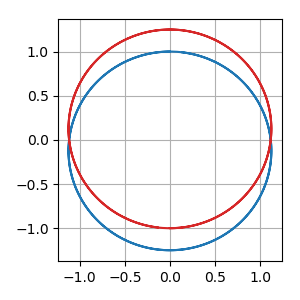}};

        %%%%
        \node [image,right=of frame11] (frame12)
        {\includegraphics[width=1.0\linewidth]{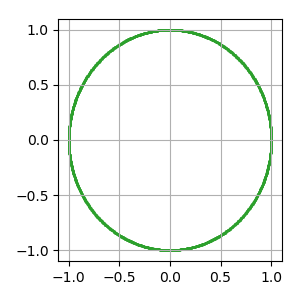}};
        
        \node [image,below=0.5cm of frame12] (frame22)
        {\includegraphics[width=1.0\linewidth]{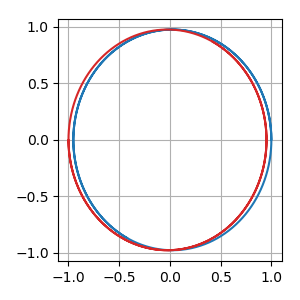}};
        
        \node [image,below=0.5cm of frame22] (frame32)
        {\includegraphics[width=1.0\linewidth]{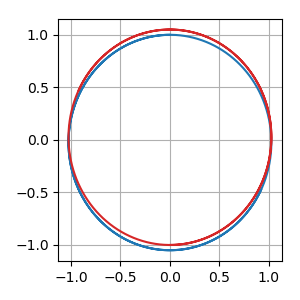}};
        
        %%%%
        \node [image,right=of frame12] (frame13)
        {\includegraphics[width=1.0\linewidth]{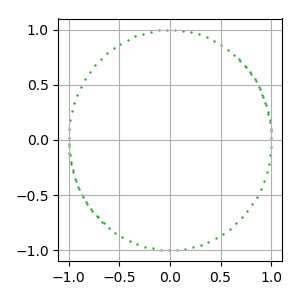}};
        
        \node [image,below=0.5cm of frame13] (frame23)
        {\includegraphics[width=1.0\linewidth]{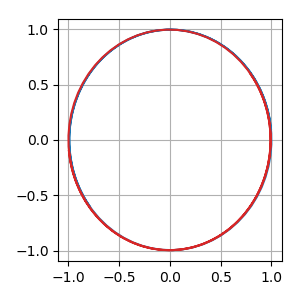}};
        
        \node [image,below=0.5cm of frame23] (frame33)
        {\includegraphics[width=1.0\linewidth]{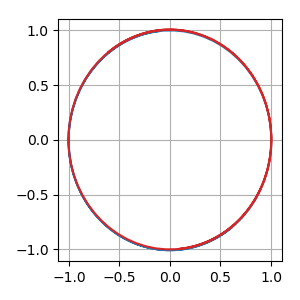}};
        
        %%%%
        \node [image,right=of frame13] (frame14)
        {\includegraphics[width=1.0\linewidth]{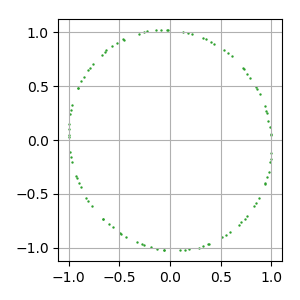}};
        
        \node [image,below=0.5cm of frame14] (frame24)
        {\includegraphics[width=1.0\linewidth]{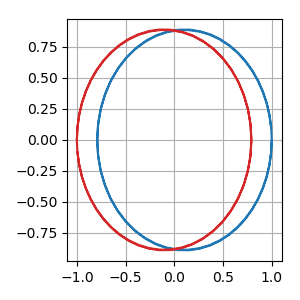}};
        
        \node [image,below=0.5cm of frame24] (frame34)
        {\includegraphics[width=1.0\linewidth]{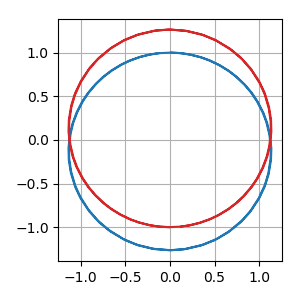}};
        
        %%%%
        \node [image,right=of frame14] (frame15)
        {\includegraphics[width=1.0\linewidth]{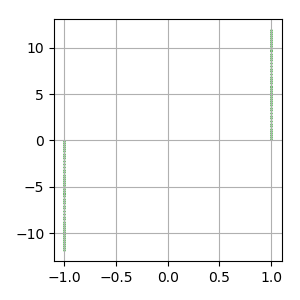}};
        
        \node [image,below=0.5cm of frame15] (frame25)
        {\includegraphics[width=1.0\linewidth]{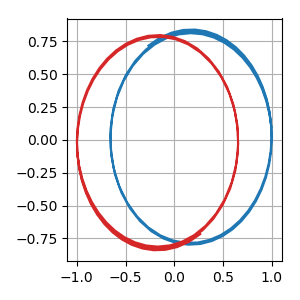}};
        
        \node [image,below=0.5cm of frame25] (frame35)
        {\includegraphics[width=1.0\linewidth]{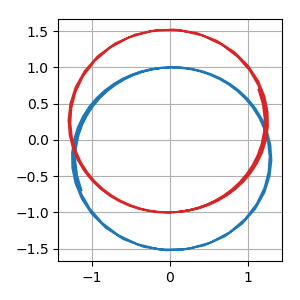}};
        
        %%%%
        \node [left=0.1cm of frame11, rotate=90,scale=0.7,xshift=0.3cm] {$p$};
        \node [below=-0.1cm of frame11 ,scale=0.7,xshift=0.3cm] {$q$};
        \node [below=-0.1cm of frame12 ,scale=0.7,xshift=0.3cm] {$q$};
        \node [below=-0.1cm of frame13 ,scale=0.7,xshift=0.3cm] {$q$};
        \node [below=-0.1cm of frame14 ,scale=0.7,xshift=0.3cm] {$q$};
        \node [below=-0.1cm of frame15 ,scale=0.7,xshift=0.3cm] {$q$};
        
        \node [left=0.1cm of frame21, rotate=90,scale=0.7,xshift=0.4cm] {$Q_y$};
        \node [below=-0.1cm of frame21 ,scale=0.7,xshift=0.3cm] {$Q_x$};
        \node [below=-0.1cm of frame22 ,scale=0.7,xshift=0.3cm] {$Q_x$};
        \node [below=-0.1cm of frame23 ,scale=0.7,xshift=0.3cm] {$Q_x$};
        \node [below=-0.1cm of frame24 ,scale=0.7,xshift=0.3cm] {$Q_x$};
        \node [below=-0.1cm of frame25 ,scale=0.7,xshift=0.3cm] {$Q_x$};
        
        \node [left=0.1cm of frame31, rotate=90,scale=0.7,xshift=0.4cm] {$P_y$};
        \node [below=-0.1cm of frame31 ,scale=0.7,xshift=0.3cm] {$P_x$};
        \node [below=-0.1cm of frame32 ,scale=0.7,xshift=0.3cm] {$P_x$};
        \node [below=-0.1cm of frame33 ,scale=0.7,xshift=0.3cm] {$P_x$};
        \node [below=-0.1cm of frame34 ,scale=0.7,xshift=0.3cm] {$P_x$};
        \node [below=-0.1cm of frame35 ,scale=0.7,xshift=0.3cm] {$P_x$};
        
        \node[rectangle,above=of frame11,xshift=0.3cm] (label1) {(a)};
        \node[rectangle,above=of frame12,xshift=0.3cm] (label2) {(b)};
        \node[rectangle,above=of frame13,xshift=0.3cm] (label3) {(c)};
        \node[rectangle,above=of frame14,xshift=0.3cm] (label4) {(d)};
        \node[rectangle,above=of frame15,xshift=0.3cm] (label5) {(e)};
    \end{tikzpicture}

    % \begin{tikzpicture}[image/.style = {text width=0.19\textwidth,inner sep=0pt, outer sep=0pt},node distance = 0mm and 0mm] 
    %     \node [image] (frame1)
    %     {\includegraphics[width=1.0\linewidth]{figs/ex2/midpoint_h=1e-1_T=10e3_eps=1e-3.eps}};
        
    %     \node [image,right=of frame1] (frame2) 
    %     {\includegraphics[width=1.0\linewidth]{figs/ex2/midpoint_h=4_T=10e3_eps=1e-3.eps}};
        
    %     \node [image,right=of frame2] (frame3) 
    %     {\includegraphics[width=1.0\linewidth]{figs/ex2/midpoint_h=100_T=10e3_eps=1e-3.eps}};
    
    %     \node [image,right=of frame3] (frame4) 
    %     {\includegraphics[width=1.0\linewidth]{figs/ex2/fastslow_h=100_T=10e3_eps=1e-3_irrational.eps}};
        
    %     \node [image,right=of frame4] (frame5) 
    %     {\includegraphics[width=1.0\linewidth]{figs/ex2/fastslow_h=100_T=10e3_eps=1e-3_resonant_pi.eps}};
        
    %     \node[rectangle,above=of frame1] (label1) {(a)};
    %     \node[rectangle,above=of frame2] (label2) {(b)};
    %     \node[rectangle,above=of frame3] (label3) {(c)};
    %     \node[rectangle,above=of frame4] (label4) {(d)};
    %     \node[rectangle,above=of frame5] (label5) {(e)};
    % \end{tikzpicture}
    
    \caption{Numerical solutions of the hidden-variable Nowtonian gravity example for $\epsilon=0.001$ and maximum integration time $T=10000$. The green trajectories on the top row denote the fast variable pair $(q,p)$. The blue and red trajectories on the mid row refer to the positions $(Q_x,Q_y)$ of the particles one and two in Cartesian coordinates, and to their respective momentums $(P_x,P_y)$ on the bottom row. The columns (a), (b), and (c) correspond to implicit midpoint scheme with time steps $h=0.1$, $h=4.0$, and $h=100.0$, respectively. The columns (d) and (e) correspond to the fast-slow scheme with a time step of $h=100.0$ and the angle variable being (d) non-resonant $\theta_0=2.0$ and (e) resonant $\theta_0=\pi$. The column (a) can be considered as the reference solution which the non-resonant fast-slow scheme reproduces quite well in the column (d).}
    \label{fig:illustration_example2}
\end{figure}

The behaviour of the scheme \eqref{eq:ex2_fast_slow_small_p}--\eqref{Bim_toy2} is illustrated in Figure \ref{fig:illustration_example2} together with the numerical solution from the well-known implicit midpoint scheme which is symplectic for canonical Hamiltonian systems and generally considered a good scheme for stiff problems. For both integrators, we set the system parameter to $\epsilon=0.001$. In Figure \ref{fig:illustration_example2}, the columns (a), (b), and (c) correspond to implicit midpoint scheme with time steps $h=0.1$, $h=4.0$, and $h=100.0$, respectively. The columns (d) and (e) correspond to the fast-slow scheme with a time step of $h=100.0$ and the angle variable being (d) non-resonant $\theta_0=2.0$ and (e) resonant $\theta_0=\pi$. The column (a) can be considered as the reference solution, which the non-resonant fast-slow integrator in column (d) closely matches. 

Certain peculiar behaviour is evident in the columns (b), (c), and (e). By increasing the time step and ``stepping over'' the fastest stiff time scale, the implicit midpoint method decouples the fast and slow variables but in an incorrect manner: through the columns (a), (b), and (c), the blue and red orbits gradually transit into co-centric circles, indicating the dependence of the adiabatic invariant $\mu_0$ on the step size $h$. Explanation for this behaviour is rooted in the asymptotic behaviour of the scheme, which is made transparent by writing the implicit midpoint scheme in the form
\begin{align}
    \overline{q}+q  &=-2\frac{\overline{p}-p}{h} -\epsilon\,4\,\frac{\overline{q}+q}{2}\,W\left(\frac{\bm{Q}+\bm{\overline{Q}}}{2}\right),\\
    \overline{p}+p & = 2\frac{\overline{q} - q}{h},\\
    \overline{\bm{P}} - \bm{P} & = -h\epsilon\, V^\prime\left(\frac{\bm{Q}+\overline{\bm{Q}}}{2}\right) - h\epsilon\,\frac{(\overline{q}+q)^2}{4}\,W^\prime\left(\frac{\bm{Q}+\overline{\bm{Q}}}{2}\right),\\
    \overline{\bm{Q}} - \bm{Q} & = h\epsilon \frac{\bm{P}+\overline{\bm{P}}}{2}.
\end{align}
When $\epsilon$ is small and the time step $h$ becomes large, the limiting behaviour of the fast variables is constant flipping of their signs. Most importantly, the term $\overline{q}+q$ becomes successively smaller with increasing $h$. Consequently, the force from the coupling potential $W$ effectively drops out from the equation for the slow variables $\bm{P}$, resulting in nearly co-centric slow orbits.

Also in column (e), some strange behaviour occurs. Via \eqref{eq:ex2_fast_slow_small_p} and \eqref{eq:ex2_fast_slow_small_q}, the resonant value of $\theta_0=\pi$ results in the following map for the fast variables
\begin{align}
     \overline{p} & = -p +\hbar\,2q\,W\left(\frac{\bm{Q}+\overline{\bm{Q}}}{2}\right),\\
     \overline{q} & = -q,
\end{align}
displaying a constant flipping of $q$, which translates to amplifying flipping in $p$. This behaviour is straightforward to verify for the initial condition $(q,p)=(1.0, 0.0)$. The behaviour of the slow variables in the column (e) in Figure \ref{fig:illustration_example2} is off from the reference solution and the non-resonant case but, because the dependence of \eqref{Bim_toy1} is only on $q$ and not on the midpoint $(\overline{q}+q)/2$, the solution still exhibits some effect from the coupling potential. Specifically, as the $q^2$ remains constant and does not average to $(q^2+p^2)/2$, the effect is actually double that of the one in the reference solution, resulting in the slow orbits being further apart from each other than what they should be.

This example serves to illustrate that even the decorated implicit midpoint scheme is not guaranteed to result in correct asymptotic behaviour and that care should be taken in trying to ``step over'' the stiff time scales. On the other hand, the example also illustrates that an integrator with the correct asymptotic behaviour may be constructed, although care is needed in choosing the saturation value for phase angle for the limit of the nearly periodic map.

\subsection{Reduced guiding-center motion\label{reduced_gc}}
We now apply the general theory developed in Section \ref{sec:symplectic_lorentz} to motion of a charged particle in a strong magnetic field of the special form $\bm{B}(x,y,z) = B(x,y)\,\bm{e}_z$, where $(x,y,z)$ denotes the usual Cartesian coordinates on $\mathbb{R}^3$ and $B$ is a positive function. 
Let $q=(x,y) \in Q=\mathbb{R}^2$ and introduce a symplectic form $\omega$ on $Q$ using the formula 
\begin{align*}
    \omega = -\mathbf{d}\alpha=-B(x,y)\,dx\wedge dy,\quad \alpha = A_x(x,y)\,dx + A_y(x,y)\,dy.
\end{align*}
Here, the components of the $1$-form $\alpha$ may be interpreted as the physicist's vector potential for $B$. Also define the Hamiltonian function $H:Q\rightarrow\mathbb{R}$ according to
\begin{align*}
    H(q) = \mu B(x,y),
\end{align*}
where $\mu$ is a positive constant parameter. The corresponding Hamiltonian vector field is given by \[X_H=\,\mathcal{R}_{\pi/2}\,\mu\nabla\ln B,\] where $\mathcal{R}_{\pi/2}$ is the rotation matrix $\mathcal{R}_{\theta}$ evaluated at $\pi/2$,
\begin{align*}
    \mathcal{R}_{\theta}=\begin{pmatrix}
      \cos\theta & -\sin\theta\\
      \sin\theta & \cos\theta\\
    \end{pmatrix}.
\end{align*}
Physically, this Hamiltonian vector field describes the motion of a charged particle's guiding center \cite{Northrop_1963} $(x,y)$. The parameter $\mu$ is the magnetic moment, and $X_H$ is also known as the $\nabla B$-drift velocity. We remark that readers familiar with the Hamiltonian formulation of guiding center motion~\cite{Li1981, CaBr2009} may be used to seeing these equations derived from the Lagrangian 
$L:TQ\rightarrow \mathbb{R}$ given by
\begin{align*}
    L(q,\dot{q})=\alpha_q(\dot{q})-H(q).
\end{align*}
We also remark that in this formulation of guiding center dynamics we have used translation invariance along $z$ to eliminate the (constant) velocity along the magnetic field and the corresponding ignorable coordinate $z$.

% and consider a Lagrangian $L:TQ\rightarrow \mathbb{R}$ given by
% \begin{align*}
%     L(q,\dot{q})=\alpha_q(\dot{q})-H(q),
% \end{align*}
% where $\alpha=A_x(x,y)dx+A_y(x,y)dy$ and $H(q)=\mu B(x,y)$. The symplectic form is given by $\omega=-\mathbf{d}\alpha=-B(x,y)\,dx\wedge dy$, and the Hamiltonian vector field by $X_H=\,\mathcal{R}_{\pi/2}\,\mu\nabla\ln B$, with $\mathcal{R}_{\pi/2}$ the rotation matrix $R_{\theta}$ evaluated at $\pi/2$
% \begin{align*}
%     R_{\theta}=\begin{pmatrix}
%       \cos\theta & -\sin\theta\\
%       \sin\theta & \cos\theta\\
%     \end{pmatrix}.
% \end{align*}
% This Lagrangian is convenient for it reproduces exactly the non-trivial contribution to the guiding-center motion of a charged particle in a magnetic field $\bm{B}=B(x,y)\bm{e}_z$: the parallel velocity coordinate remains constant, $\dot{z}=0$, and the drift motion is exactly recovered by $X_H$.

In order to construct the symplectic Lorentz map for this system, we begin by observing that the complex structure
\[
\mathbb{J}\begin{pmatrix}
\dot{x}\\
\dot{y}
\end{pmatrix} = \begin{pmatrix}
0 & 1\\
-1 & 0
\end{pmatrix}\begin{pmatrix}
\dot{x}\\
\dot{y}
\end{pmatrix}
\] 
is compatible with $\omega$ since
\begin{align*}
    \omega(\dot{q}_1,\mathbb{J}\dot{q}_2) = B(x,y)\,\dot{q}_1\cdot \dot{q}_2,
\end{align*}
where $\cdot$ denotes the usual inner product on $\mathbb{R}^2$. We may therefore use the metric $g_q(v,w)=B(x,y)v\cdot w$ to build a new Hamiltonian system on $TQ$, compatible with the Lorentz-embedding idea from \cite{Burby_Hirvijoki_2021_JMP}, 
\begin{align*}
    \Omega_{\epsilon}^\ast=\pi^\ast\omega-\epsilon\,\mathbf{d}(g_q(v,dq)), \qquad H_\epsilon^\ast(q,v)=\frac{1}{2}\epsilon^2\,g_q(v,v)+\epsilon \,\tau H(q),
\end{align*}
where we have also introduced a constant factor $\tau$ for scaling time of the original system. This scaling will help us later to assign the fastest motion in the guiding center system to occur at order one and help in nonlinear solve of the coordinate update map in our numerical example. Essentially, the true time of the original system now evolves at rate that is $\tau$ times the rate of the embedded system.
The equations of motion for this larger system are given by
\begin{align*}
    \dot{q}&=\epsilon v,\\
    \dot{v}&=-(1+\epsilon v\cdot \mathcal{R}_{\pi/2}\cdot\nabla\ln B)\,\mathcal{R}_{\pi/2}v-\left(\tau\mu+\epsilon\frac{1}{2}|v|^2\right)\nabla\ln B.
\end{align*}
%This is a nearly periodic system with unit frequency and a stable slow manifold agreeing to the lowest order with the reduced guiding-center motion. 

Proceeding now with the general construction of the symplectic Lorentz map, we introduce a Type I generating function
\begin{align*}
    S(q,\overline{q}) &= \int_q^{\overline{q}}\alpha + \Sigma(q/2+\overline{q}/2,\overline{q}-q)  ,
\end{align*}
where $\Sigma:TQ\rightarrow \mathbb{R}$ is given by
\begin{align*}
\Sigma(\eta,\xi) &= -\hbar \mu B(\eta) + \hbar^2\,B(\eta)\,X_H(\eta)\cdot\xi \nonumber\\
    & \qquad- \frac{1}{4}\left(\frac{\sin\theta_0}{1-\cos\theta_0}\right)B(\eta)\,(\xi-\hbar X_H(\eta))\cdot(\xi-\hbar X_H(\eta)).
\end{align*}
Note we are using the symbol $\eta$ instead of $x$, in contrast to Section \ref{sec:symplectic_lorentz}, in order to avoid confusion with the standard Cartesian coordinate system. The term involving derivatives of $\omega$, present in the previous section, vanishes identically for $X_H\cdot\nabla B=0$.
The symplectic Lorentz map then provides the equations
\begin{align*}
    &\hbar^2B(\overline{q})\overline{v}=-\int_0^1\lambda B(q+\lambda\xi)d\lambda\,\mathcal{R}_{\pi/2}\xi+\frac{1}{2}\partial_\eta\Sigma(\eta,\xi)+\partial_\xi\Sigma(\eta,\xi),\\
    &\hbar^2B(q)v=\int_0^1(1-\lambda) B(q+\lambda\xi)d\lambda\,\mathcal{R}_{\pi/2}\xi-\frac{1}{2}\partial_\eta\Sigma(\eta,\xi)+\partial_\xi\Sigma(\eta,\xi),
\end{align*}
where the derivatives are
\begin{align*}
    \partial_\eta\Sigma(\eta,\xi)&=-\hbar \mu \nabla B(\eta) + \hbar^2\,\nabla B(\eta)\,X_H(\eta)\cdot\xi
    + \hbar^2\,B(\eta)\,\nabla X_H(\eta)\cdot\xi \nonumber\\
    &\qquad - \frac{1}{4}\left(\frac{\sin\theta_0}{1-\cos\theta_0}\right)\nabla B(\eta)\,(\xi-\hbar X_H(\eta))\cdot(\xi-\hbar X_H(\eta))\nonumber\\
    &\qquad+\frac{1}{2}\left(\frac{\sin\theta_0}{1-\cos\theta_0}\right)B(\eta)\,\hbar \nabla X_H(\eta)\cdot(\xi-\hbar X_H(\eta)),
    \\
    \partial_\xi\Sigma(\eta,\xi) &= \hbar^2\,B(\eta)\,X_H(\eta) - \frac{1}{2}\left(\frac{\sin\theta_0}{1-\cos\theta_0}\right)B(\eta)\,(\xi-\hbar X_H(\eta)),
\end{align*}
and everything is understood to be evaluated at $(\eta,\xi)=((\overline{q}+q)/2,\overline{q}-q)$. 

Next, we rearrange the implicit equation for $\overline{q}$ into
% \begin{align*}
%     & \left[\frac{1}{2}\left(\frac{\sin\theta_0}{1-\cos\theta_0}\right)B(\eta)\,1-\int_0^1(1-\lambda) B(q+\lambda\xi)d\lambda\,\mathcal{R}_{\pi/2}\right]\left[\xi-\hbar X_H(\eta)\right]\\
%     &=\hbar \int_0^1(1-\lambda) B(q+\lambda\xi)d\lambda\,\mathcal{R}_{\pi/2}X_H(\eta)\\
%     &+\hbar^2\,B(\eta)\,X_H(\eta)-\hbar^2B(q)v
%     \\
%     &-\frac{1}{2}\partial_\eta\Sigma(\eta,\xi)
% \end{align*}
% The alternative route is 
\begin{align*}
    & \left[\frac{1}{2}\left(\frac{\sin\theta_0}{1-\cos\theta_0}\right)B(\eta)\,1-\int_0^1(1-\lambda) B(q+\lambda\xi)d\lambda\,\mathcal{R}_{\pi/2}\right]\xi\\
    &=\left[\frac{1}{2}\left(\frac{\sin\theta_0}{1-\cos\theta_0}\right)B(\eta)\,1-\frac{1}{2}B(\eta)\mathcal{R}_{\pi/2}\right]\,\hbar X_H(\eta)\\
    &\qquad
    +\hbar^2\,B(\eta)\,X_H(\eta)-\hbar^2B(q)v
    %\\&\qquad
    -\frac{1}{2} \hbar^2\,\nabla B(\eta)\,X_H(\eta)\cdot\xi
    %\\&\qquad
    -\frac{1}{2} \hbar^2\,B(\eta)\,\nabla X_H(\eta)\cdot\xi
    \\&\qquad
    + \frac{1}{8}\left(\frac{\sin\theta_0}{1-\cos\theta_0}\right)\nabla B(\eta)\,(\xi-\hbar X_H(\eta))\cdot(\xi-\hbar X_H(\eta))
    \\&\qquad
    -\frac{1}{4}\left(\frac{\sin\theta_0}{1-\cos\theta_0}\right)B(\eta)\,\hbar \nabla X_H(\eta)\cdot(\xi-\hbar X_H(\eta)),
\end{align*}
which can be iterated for $\eta$ and $\xi$. After that, one solves for $\overline{v}$ by evaluating, for example, the expression
\begin{align*}
    \hbar^2B(\overline{q})\overline{v}+\hbar^2B(q)v&=-\int_0^1\lambda B(q+\lambda\xi)d\lambda\,\mathcal{R}_{\pi/2}\xi+\int_0^1(1-\lambda) B(q+\lambda\xi)d\lambda\,\mathcal{R}_{\pi/2}\xi
    \\
    &\qquad+2\hbar^2\,B(\eta)\,X_H(\eta) - \left(\frac{\sin\theta_0}{1-\cos\theta_0}\right)B(\eta)\,(\xi-\hbar X_H(\eta)).
\end{align*}

% {\color{red}Here is an alternative nonlinear solve strategy. -JB} First write the $q$  and $\overline{q}$ components of the generating function relations as
% \begin{align*}
%     \hbar^2\,B(q)\,v = -\left(\int_0^1[1-\lambda]\,B_\lambda\,d\lambda\right)\,\mathbb{J}\xi 
% \end{align*}

\begin{figure}[!htb]
    \centering
    \begin{tikzpicture}
    \node[above right,inner sep = 0] (image) at (0,0) {
        \includegraphics[width=0.5\linewidth]{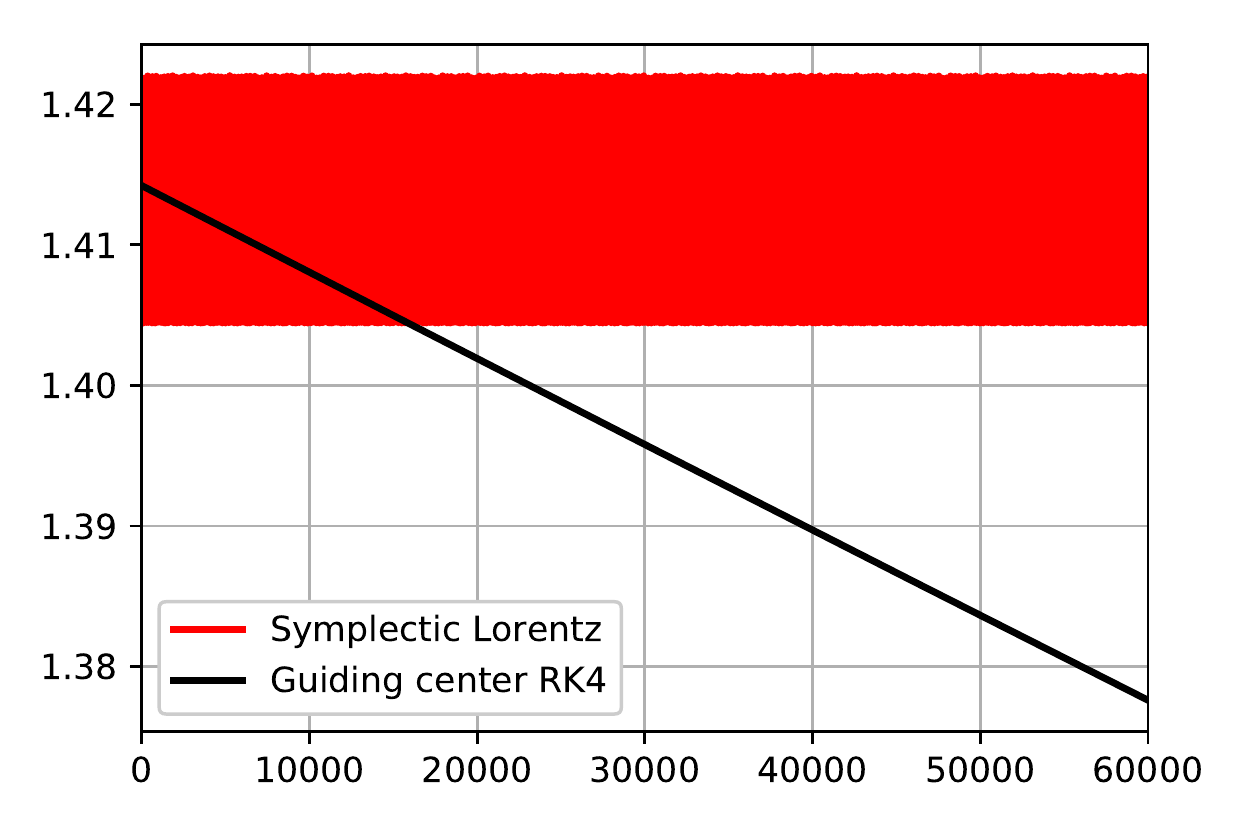}
    };
    \begin{scope}[
        x={($1*(image.south east)$)},
        y={($1*(image.north west)$)}]
        \node[rectangle,rotate=90] at (-0.03,0.5) {$|q|$};
        \node[rectangle] at (0.5,-0.03) {\# steps};
    \end{scope}
    \end{tikzpicture}
    \caption{Comparison of the guiding-center RK4 integrator and the symplectic Lorentz map in the simple magnetic field case. The orbit radius $|q|$ of the RK4 integrator deteriorates clearly while the symplectic Lorentz map manages to retain the oscillations in the radius within limits stable limits.}
    \label{fig:RK4_vs_symplectic_lorentz_simple}
\end{figure}

Next, we perform some numerical tests. First, we choose a magnetic field
\begin{align*}
    B=B_0(1+\alpha |q|^2),
\end{align*}
where $\alpha$ introduces a small perturbation to the otherwise constant magnetic field. For the original system $\dot{q}=X_H(q)$, this field results in circular orbits for $q$. We then investigate the solutions of the symplectic Lorentz map and compare them with the classic RK4 integrator applied to the original system. For the scaling of time, we choose $\tau=\alpha^{-1}$. Choosing an initial point $q=(1,1)$, parameters $B_0=1$, $\mu=1.0$, $\alpha=0.001$, and initializing with $v=X_H(q)$, we run the simulation for 60'000 steps of size $\hbar=0.1$. This is enough to demonstrate the deterioration of the RK4 method while the symplectic Lorentz map preserves the orbit in place, as seen in Figure \ref{fig:RK4_vs_symplectic_lorentz_simple}.

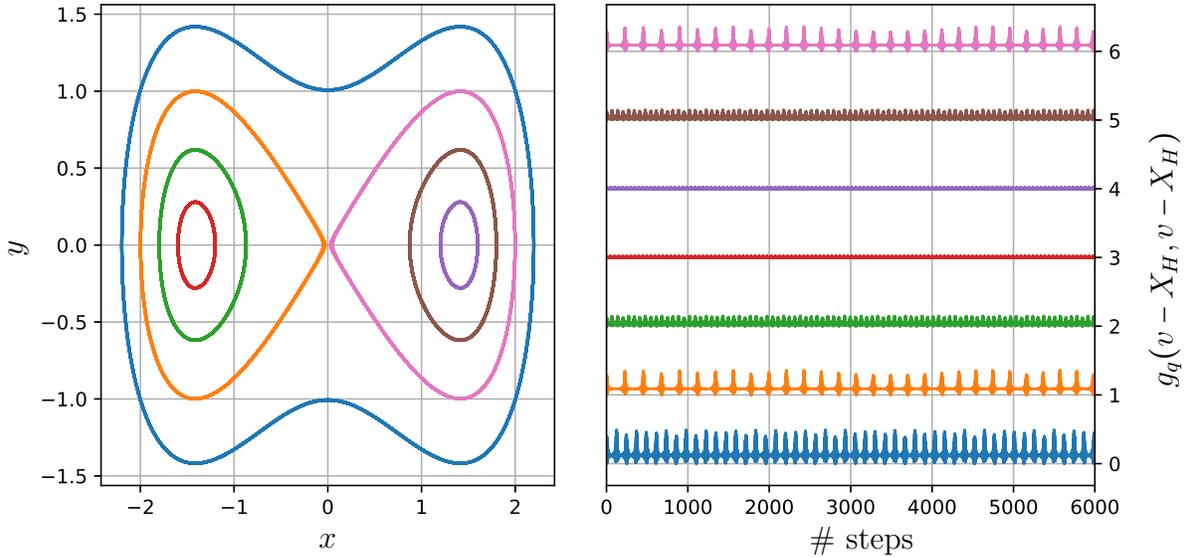
\begin{figure}[!htb]
    \centering
    \begin{tikzpicture}
    \node[above right,inner sep = 0] (image) at (0,0) {
    \begin{tikzpicture}[image/.style = {text width=0.45\textwidth,inner sep=0pt, outer sep=0pt},node distance = 0mm and 0mm] 
        \node [image] (frame1)
        {\includegraphics[width=\textwidth]{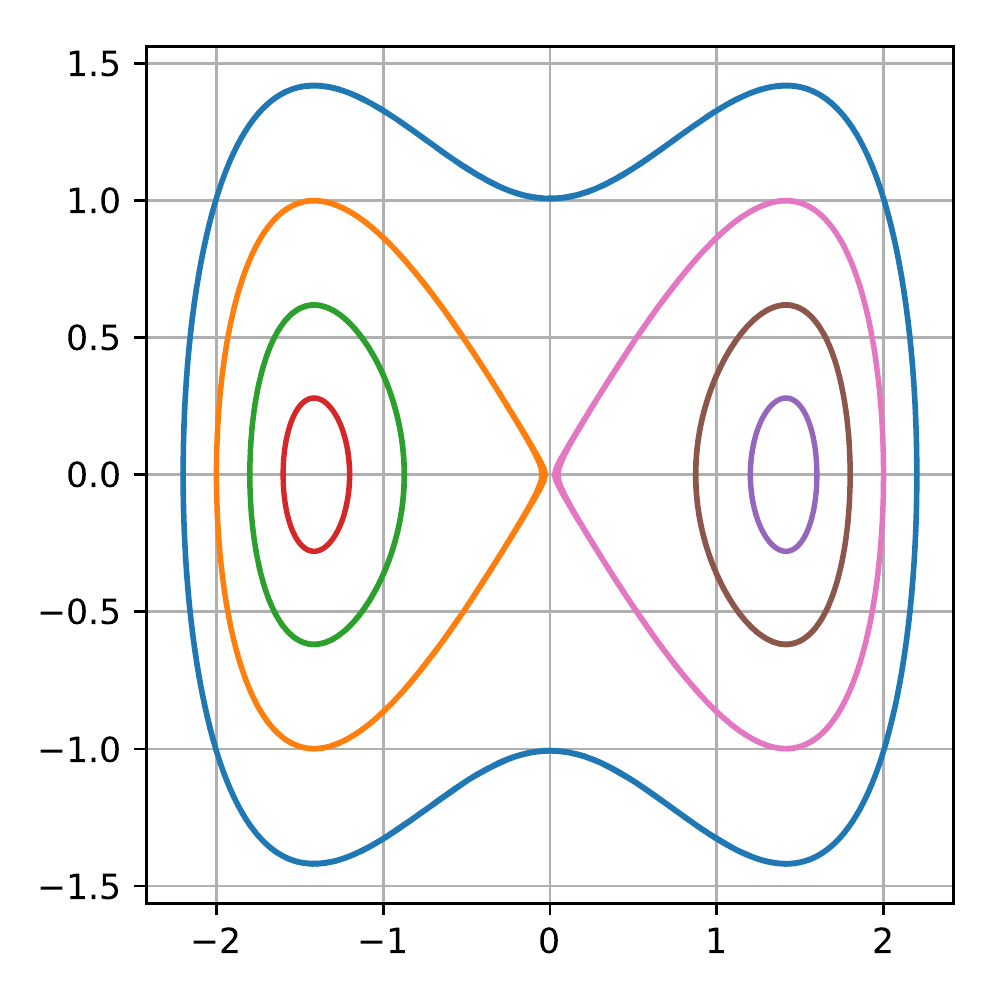}};
        \node [image,right=of frame1] (frame2) 
        {\includegraphics[width=\textwidth]{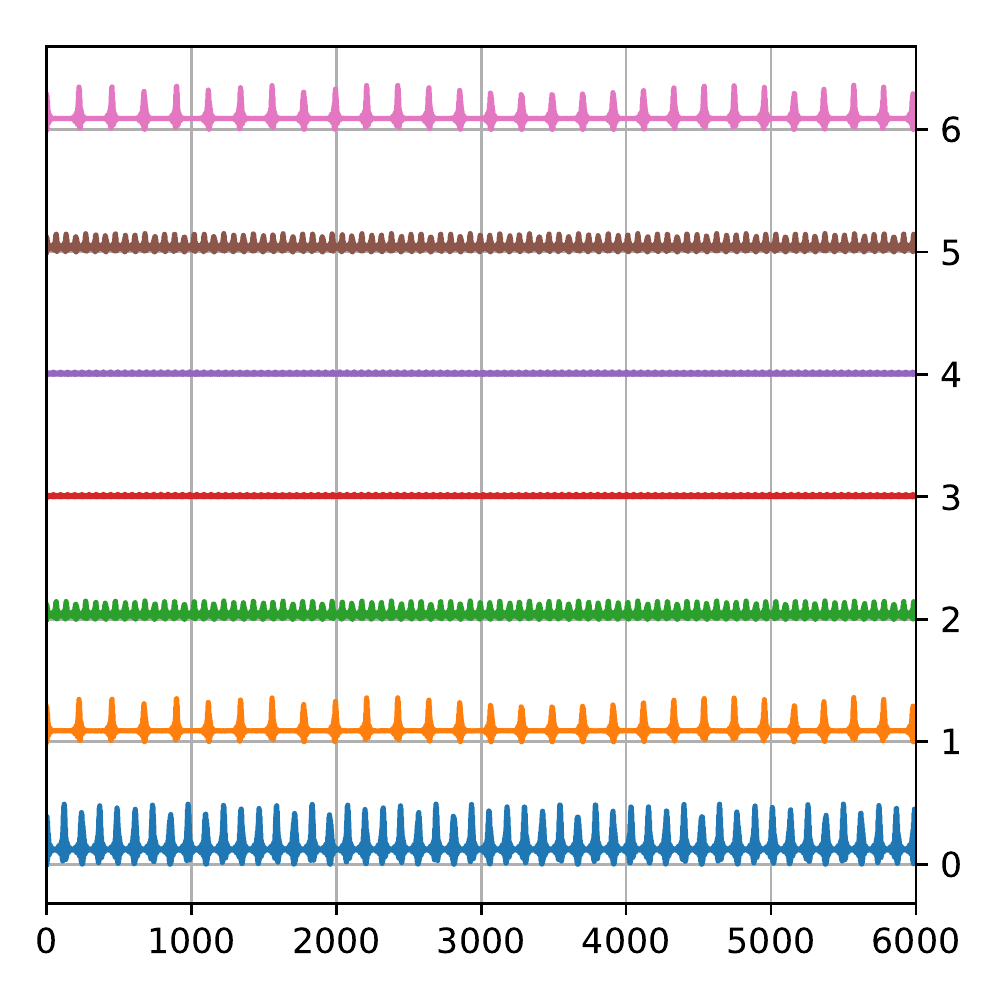}};
    \end{tikzpicture}
    };
    \begin{scope}[
        x={($1*(image.south east)$)},
        y={($1*(image.north west)$)}]
        %\draw[red,step=0.1] (image.south west) grid (image.north east);
        \node[rectangle] at (0.275,0.) {$x$};
        \node[rectangle,rotate=90] at (0.0,0.52) {$y$};
        \node[rectangle] at (0.75, 0.0) {\# steps};
        \node[rectangle,rotate=90] at (1.02,0.5) {$g_q(v-X_H,v-X_H)$};
    \end{scope}  
    \end{tikzpicture}
    \caption{Phase-space orbits (left) and the adiabatic invariant (right). For the ``figure-eight'' magnetic field. The values of the adiabatic invariant have been shifted by the integers $\{0,...,6\}$ for illustrative purposes. A lesser number of steps (6000) have been chosen to illustrate the non-trivial beating structure of the adiabatic invariant.}
    \label{fig:reduced_figure_eight}
\end{figure}

Next we consider the magnetic field
\begin{align*}
    B(x,y)=2+y^2-x^2+\frac{1}{4}x^4,
\end{align*}
whose level sets have a ``figure-eight'' structure. By energy conservation, the guiding-center orbits should reflect this pattern. Choosing a time step of $\hbar=0.05$ and $\tau=1.0$, we run the symplectic Lorentz map for 6'000 steps and illustrate both the orbits and the evolution of the postulated adiabatic invariant $g_q(v-X_H,v-X_H)$ in Figure~\ref{fig:reduced_figure_eight}. The orbits appear stable and well confined to their respective phase-space domains, and the adiabatic invariant remains within bounds while oscillating with a non-trivial beating structure.

\begin{figure}[!htb]
    \centering
    \begin{tikzpicture}
    \node[above right,inner sep = 0] (image) at (0,0) {
    \begin{tikzpicture}[image/.style = {text width=0.45\textwidth,inner sep=0pt, outer sep=0pt},node distance = 0mm and 0mm] 
        \node [image] (frame1)
        {\includegraphics[width=0.97\textwidth]{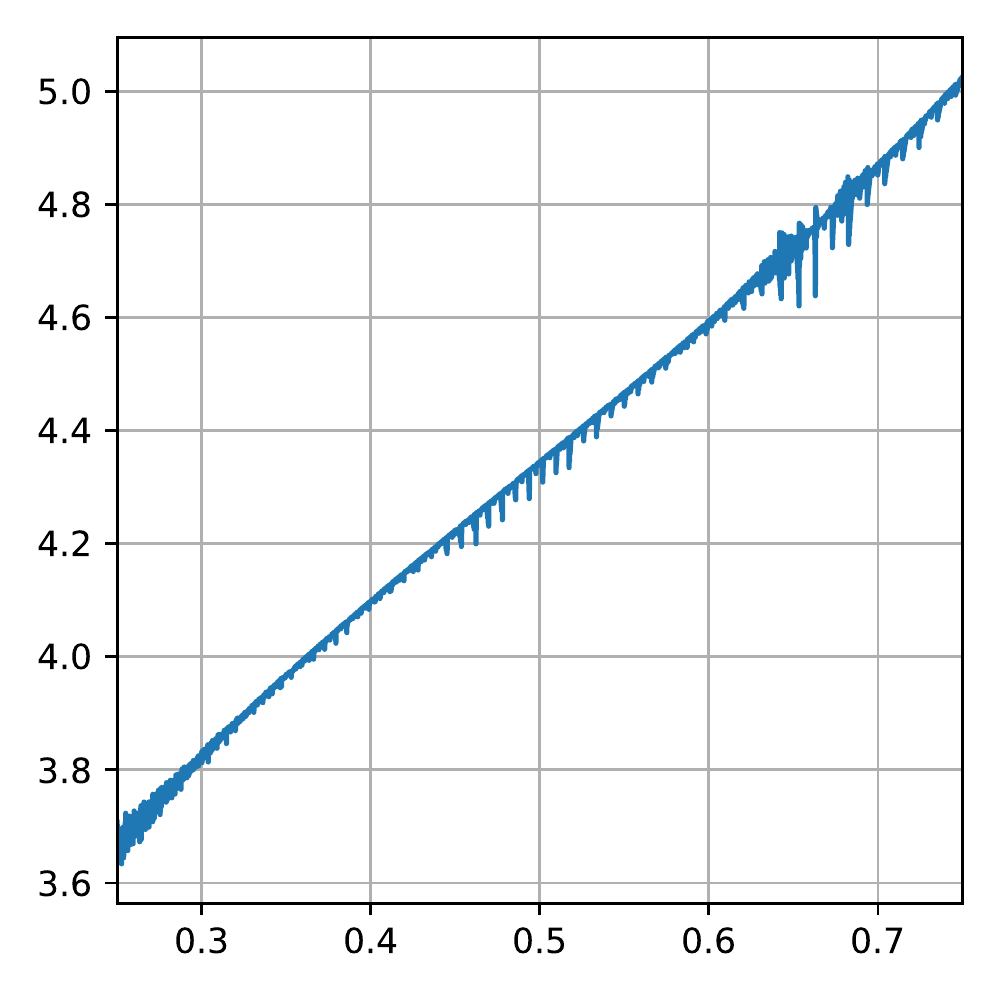}};
        \node [image,right=of frame1] (frame2) 
        {\includegraphics[width=\textwidth]{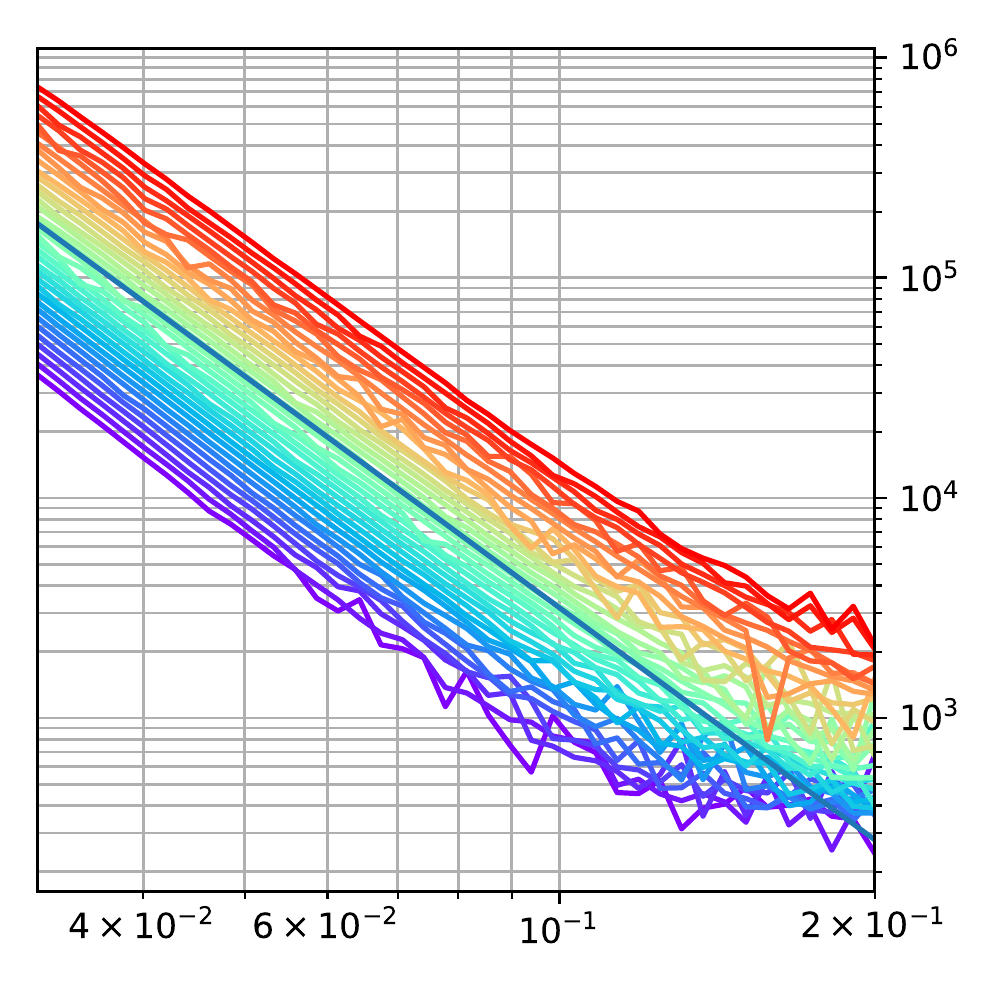}};
    \end{tikzpicture}
    };
    \begin{scope}[
        x={($1*(image.south east)$)},
        y={($1*(image.north west)$)}]
        %\draw[red,step=0.1] (image.south west) grid (image.north east);
        \node[rectangle] at (0.275,0.) {$\theta_0/\pi$};
        \node[rectangle,rotate=90] at (0.0,0.55) {$\text{log}_{10}T_\text{breakdown}$};
        \node[rectangle] at (0.75, 0.0) {$\hbar$};
        \node[rectangle,rotate=90] at (1.00,0.52) {$T_\text{breakdown}$};
    \end{scope}  
    \end{tikzpicture}

    \caption{Left panel: Breakdown time for the adiabatic invariant versus $\theta_0$. Timestep is fixed at $\hbar = 0.05565$. Initial condition is $x = 2.0$, $y = 0.0$. Right panel: Breakdown time for the adiabatic invariant versus $\hbar$. Colorscale indicates value of $\theta_0$, ranging from $\theta_0 = \pi/4$ (purple) to $\theta_0 = 3\pi/4$ (red). Initial condition is $x = 2.0$, $y = 0.0$. Central dark green line is $\hbar^{-3.5}$, for reference. Theory predicts that the breakdown time should scale like $\hbar^{-N}$ for any positive $N$ when $\hbar$ is small enough. While superpolynomial scaling of the breakdown time as a function of $\hbar$ is not apparent in these computations, it cannot be ruled out given the limited range of $\hbar$ values considered. Computing the breakdown time for appreciably smaller values of $\hbar$ rapidly becomes prohibitively expensive because the adiabatic invariant is so well conserved. }
    \label{fig:my_label}
\end{figure}
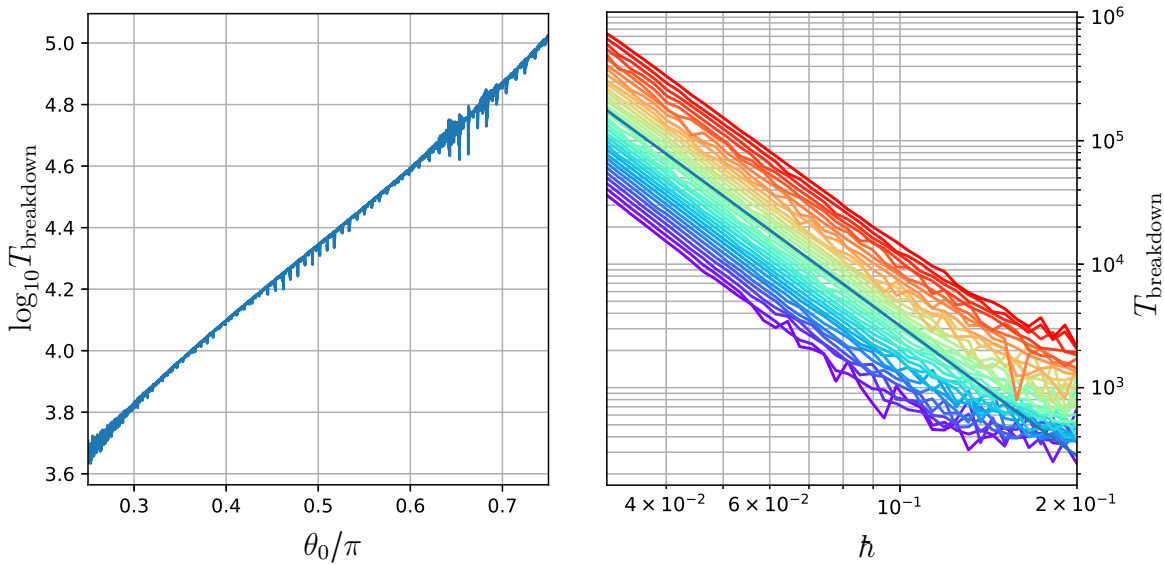

For the same magnetic field, we performed a pair of tests that probe the robustness of the discrete-time adiabatic invariant $\mu$. We introduce an empirical estimate of the breakdown time for $\mu$ conservation and compute how that estimate varies with the parameters $\theta_0$ and $\hbar$. Our estimate is based on the observation that $\mu$ typically oscillates about a time-varying mean value $\overline{\mu}$ with an approximately constant oscillation amplitude $\tilde{\mu}$.  We estimate that breakdown has occured after $n$ iterations when $\overline{\mu}(n\,\hbar)-\overline{\mu}(0)> \tilde{\mu}$. We then define the breakdown time estimate to be $T_{\text{breakdown}} = n\,\hbar$. Results from our sensitivy studies are displayed in Figure~\ref{fig:my_label}. While the general theory predicts that the breakdown time should scale as fast as $\hbar^{-N}$ for any non-negative integer $N$, the observable asymptote in $T_{\text{breakdown}}(\hbar)$ appears well-approximated by $\hbar^{-3.5}$. We presently lack understanding of the origin of the scaling exponent $-3.5$. Theory also predicts that adiabatic invariance should be less robust when $\theta_0/2\pi$ is rational. This prediction is consistent with the plot of $T_{\text{breakdown}}(\theta_0)$, which shows intermittent depressions in the breakdown time superposed on a strong upward trend as $\theta_0$ approaches $\pi$. We hypothesize that these depressions occur at small denominator rational values of $\theta_0/2\pi$ that produce nonlinear self-resonance in the integrator. As with the scaling exponent, we presently lack detailed understanding for the dramatic increase in observed breakdown time as $\theta_0$ approaches $\pi$. 

% \begin{figure}
%     \centering
%     \includegraphics[scale=1]{figs/Tb_vs_h_loglog.eps}
%     \caption{Breakdown time for the adiabatic invariant versus $\hbar$. Colorscale indicates value of $\theta_0$, ranging from $\theta_0 = \pi/4$ (purple) to $\theta_0 = 3\pi/4$ (red). Initial condition is $x = 2.0$, $y = 0.0$. Central dark green line is $\hbar^{-3.5}$, for reference. Theory predicts that the breakdown time should scale like $\hbar^{-N}$ for any positive $N$ when $\hbar$ is small enough. Therefore the sampled values of $\hbar$ are not small enough to be in asymptotic regime. Computing the breakdown time for appreciably smaller values of $\hbar$ is prohibitively expensive. }
%     \label{fig:my_label}
% \end{figure}

\section{Discussion\label{discussion_section}}
In this article we have introduced and developed the theoretical foundations of nearly-periodic maps. These maps provide a discrete-time analogue of Kruskal's \cite{Kruskal_1962} continuous-time nearly-periodic systems. The limiting dynamics of both nearly-periodic systems and nearly-periodic maps translate points along the orbits of a principal circle bundle. In the continuous-time case, each limiting trajectory ergodically samples an orbit. In discrete time, non-resonance appears as an additional requirement for ergodic sampling. As a first major application of nearly-periodic maps, we used them to construct a class of geometric integrators for Hamiltonian systems on arbitrary exact symplectic manifolds. 

Kruskal's principal interest in continuous-time nearly-periodic systems came from their relationship to the theory of adiabatic invariants. In the paper \cite{Kruskal_1962}, Kruskal showed that nearly-periodic systems necessarily admit approximate $U(1)$-symmetries. He then went on to deduce that this approximate symmetry implies the existence of an adiabatic invariant when the underlying nearly-periodic system happens to be Hamiltonian. The theory of nearly-periodic maps is satisfying in this respect since it establishes the existence of a discrete-time adiabatic invariant for nearly-periodic maps with an appropriate Hamiltonian structure. Moreover, the arguments used in the existence proof parallel those originally used by Kruskal. (See Thm. \ref{nps_noether}.)

It is useful to place the integrators developed in this article in the context of previous attempts at geometric integration of noncanonical Hamiltonian systems. Based on the observation \cite{Arnold_1989} that Hamiltonian systems on exact symplectic manifolds admit degenerate ``phase space Lagrangians" \cite{Cary_1983}, Qin \cite{Qin_2008} proposed direct application of the theory of variational integration \cite{MaWe2001} to phase space Lagrangians for noncanonical systems. While initial results looked promising, further investigations by Ellison \cite{Ellison_thesis,Ellison_2018} revealed that the most intuitive variational discretizations of phase space Lagrangians typically suffer from unphysical instabilities known as ``parasitic modes" \cite{Hairer_2006}. As noticed first in \cite{Rowley_2002}, the origin of these parasitic modes is related to a mismatch between the differing levels of degeneracy in the phase space Lagrangian and its discretization. Our integrators may be understood as modifications of those studied by Qin and Ellison that stablize the parasitic modes over very large time intervals by way of a discrete-time adiabatic invariant. This ``adiabatic stabilization" mechanism is conceptually interesting since it suppresses numerical instabilities without resorting to the addition of artificial dissipation. Also of note, adiabatic stablization differs from the stablization mechanism proposed by Ellison in \cite{Ellison_2018}, wherein the phase space Lagrangian is discretized so that it has the same level of degeneracy as its continuous-time counterpart. While Ellison's ``properly-degenerate" discretizations apply to a very limited class of non-canonical Hamiltonian systems, (see \cite{Ellison_2018} for the precise limitations) the adiabatic stablization method discussed here applies to any Hamiltonian system on an exact symplectic manifold.

In the preprint \cite{Kraus_2017_arXiv}, Kraus has developed an alternative approach to structure-preserving integration of noncanonical Hamiltonian systems based on projection methods. In contrast to our approach, this technique is designed to produce integrators that preserve the original system's symplectic form, rather than a symplectic form on a larger space. However, there is no geometric picture for why Kraus' method ought to have this property. In fact, Kraus finds that geometrically-reasonable variants of his method are not symplectic. The structure-preserving properties of our method are easier to understand in this respect, since they follow from the standard theory of mixed-variable generating functions for symplectic maps. Both techniques warrant further investigation.

As a final remark concerning relationships between the theory developed here and previous work, it is worthwhile highlighting the technique introduced by Tao in \cite{Tao_2016} for constructing explicit symplectic integrators for non-separable Hamiltonians. The latter technique applies to canonical Hamiltonian systems with general Hamiltonian $H(q,p)$. It proceeds by constructing a canonical Hamiltonian system in a space with double the dimension of the original $(q,p)$ space, and then applying splitting methods to the larger system. Much like the symplectic Lorentz system introduced in \cite{Burby_Hirvijoki_2021_JMP}, and exploited in Section \ref{sec:symplectic_lorentz}, Tao's larger system contains a copy of the original system as a normally-elliptic invariant manifold. This suggests that Tao's construction might be interpreted as an application of nearly-periodic maps. It is a curious fact, however, that Tao's error analysis suggests the oscillation frequency around the invariant manifold cannot be made to be arbitrarily large. This indicates nearly-periodic map theory is not an appropriate tool for understanding Tao's results. It would be interesting to investigate whether or not nearly-periodic map theory can be used to sharpen Tao's estimates.

More work is required to develop nearly-periodic map machinery, in both theory and practice. The following is a list of just a few of the open theoretical questions in this area.
\begin{enumerate}
    \item Non-resonant nearly-periodic maps and nearly-periodic systems admit formal $U(1)$-symmetries, and therefore formal reductions to the space of $U(1)$-orbits. Given an arbitrary continuous-time nearly-periodic system, are there systematic strategies for constructing nearly-periodic maps whose $U(1)$-reduction approximates the flow of the nearly-periodic system's $U(1)$-reduction? (We provide a simple example where this can be done in Section \ref{hidden_gravity_section}.) Such maps would providee a good approximation of the original system's dynamics ``on average."
    \item For a Hamiltonian system on an exact symplectic manifold $M$, the geometric integrators constructed in this Article comprise symplectic mappings on $TM$ that admit approximate invariant manifolds diffeomorphic to $M$. In light of this diffeomorphism, is there some sense in which our integrators possess an adiabatically-invariant symplectic form on $M$? (Note that this does not obviously follow from the symplectic property on $TM$.) 
    \item A commonly touted benefit of symplectic integration is the long-time approximate preservation of energy. Proofs of this result rely on backward error analysis. Can similar techniques be used to prove that our geometric integrators approximately preserve the original Hamiltonian system's energy, at least over large time intervals? Our initial numerical experiments suggest such a result is satisfied for as long as the discrete-time adiabatic invariant is well-conserved.
    \item Our geometric integrators enjoy local $O(\hbar^{5/2})$ accuracy. Are there extensions of these integrators with arbitrarily-high-order accuracy?
\end{enumerate}

\section*{Acknowledgements}

The work of JWB was supported by the Los Alamos National Laboratory LDRD program under Project No.~20180756PRD4. The work of EH was supported by the Academy of Finland (Grant No.~315278). Any subjective views or opinions expressed herein do not necessarily represent the views of the Academy of Finland or Aalto University. The work of ML was supported by NSF under grants DMS-1345013, DMS-1813635, by AFOSR under grant FA9550-18-1-0288, and by the DoD under grant HQ00342010023 (Newton Award for Transformative Ideas during the COVID-19 Pandemic).

\bibliographystyle{unsrt}
\bibliography{cumulative_bib_file.bib}

\providecommand{\noopsort}[1]{}\providecommand{\singleletter}[1]{#1}%
\begin{thebibliography}{10}

\bibitem{Kruskal_1962}
M.~Kruskal.
\newblock {Asymptotic theory of Hamiltonian and other systems with all
  solutions nearly periodic}.
\newblock {\em J. Math. Phys.}, 3:806, 1962.

\bibitem{Burby_Squire_2020}
J.~W. Burby and J.~Squire.
\newblock General formulas for adiabatic invariants in nearly periodic
  hamiltonian systems.
\newblock {\em J. Plasma Phys.}, 86:835860601, 2020.

\bibitem{Burby_Hirvijoki_2021_JMP}
J.~W. Burby and H.~Hirvijoki.
\newblock {Normal stability of slow manifolds in nearly-periodic Hamiltonian
  systems}.
\newblock {\em J. Math. Phys.}, 62:093506, 2021.

\bibitem{Ricketson_2020}
L.~F. Ricketson and L.~Chac\'on.
\newblock An energy-conserving and asymptotic-preserving charged-particle orbit
  implicit time integrator for arbitrary electromagnetic fields.
\newblock {\em Journal of Computational Physics}, 418:109639, Oct 2020.

\bibitem{Xiao_Qin_2021}
J.~Xiao and H.~Qin.
\newblock Slow manifolds of classical pauli particle enable
  structure-preserving geometric algorithms for guiding center dynamics.
\newblock {\em Computer Physics Communications}, 265:107981, Aug 2021.

\bibitem{Hairer_Lubich_2020}
E.~Hairer and C.~Lubich.
\newblock {Long-time analysis of a variational integrator for charged particle
  dynamics in a strong magnetic field}.
\newblock {\em Num. Math.}, 144:699--728, 2020.

\bibitem{Abraham_2008}
R.~Abraham and J.~E. Marsden.
\newblock {\em Foundations of Mechanics}.
\newblock AMS Chelsea publishing. American Mathematical Soc., 2008.

\bibitem{Marsden_Ratiu_1999}
J.~E. Marsden and T.~S. Ratiu.
\newblock {\em Introduction to Mechanics and Symmetry}, volume~17 of {\em Texts
  in Applied Mathematics}.
\newblock Springer, New York, NY, 1999.

\bibitem{Northrop_1963}
T.~G. Northrop.
\newblock {\em The Adiabatic Motion of Charged Particles}.
\newblock Interscience tracts on physics and astronomy. Interscience
  Publishers, 1963.

\bibitem{Littlejohn_1981}
R.~G. Littlejohn.
\newblock Hamiltonian formulation of guiding center motion.
\newblock {\em Phys.\ Fluids}, 24:1730, 1981.

\bibitem{Littlejohn_1983}
R.~G. Littlejohn.
\newblock Variational principles of guiding centre motion.
\newblock {\em J. Plasma Phys.}, 29:111, 1983.

\bibitem{Littlejohn_1979}
R.~G. Littlejohn.
\newblock A guiding center hamiltonian: A new approach.
\newblock {\em J. Math. Phys.}, 20:2445, 1979.

\bibitem{Lochak_1993}
P.~Lochak.
\newblock Hamiltonian perturbation theory: periodic orbits, resonances and
  intermittency.
\newblock {\em Nonlinearity}, 6:885--904, 1993.

\bibitem{Silva_book_2008}
A.~C. da~Silva.
\newblock {\em Lectures on Symplectic Geometry}.
\newblock Lecture Notes in Mathematics. Springer, 2008.

\bibitem{Burby_Klotz_2020}
J.~W. Burby and T.~J. Klotz.
\newblock Slow manifold reduction for plasma science.
\newblock {\em Comm. Nonlin. Sci. Numer. Simul.}, 89:105289, 2020.

\bibitem{Kraus_2017_arXiv}
M.~Kraus.
\newblock Projected variational integrators for degenerate {L}agrangian
  systems, 2017.

\bibitem{MaWe2001}
{J.E.} Marsden and M.~West.
\newblock Discrete mechanics and variational integrators.
\newblock {\em Acta Numer.}, 10:317--514, 2001.

\bibitem{LeZh2009}
M.~Leok and J.~Zhang.
\newblock Discrete {H}amiltonian variational integrators.
\newblock {\em IMA J. Numer. Anal.}, 31(4):1497--1532, 2011.

\bibitem{MuOr2004}
S.~M{\"u}ller and M.~Ortiz.
\newblock On the $\gamma$-convergence of discrete dynamics and variational
  integrators.
\newblock {\em Journal of Nonlinear Science}, 14(3):279--296, 2004.

\bibitem{HaLe2015}
J.~Hall and M.~Leok.
\newblock Spectral variational integrators.
\newblock {\em Numerische Mathematik}, 130(4):681--740, 2015.

\bibitem{IsNo2004}
A.~Iserles and S.~P. N{\o}rsett.
\newblock On quadrature methods for highly oscillatory integrals and their
  implementation.
\newblock {\em BIT Numerical Mathematics}, 44(4):755--772, 2004.

\bibitem{StGr2009}
A.~Stern and E.~Grinspun.
\newblock Implicit-explicit variational integration of highly oscillatory
  problems.
\newblock {\em Multiscale Model. Simul.}, 7(4):1779--1794, 2009.

\bibitem{Li1981}
R.~G. Littlejohn.
\newblock Hamiltonian formulation of guiding center motion.
\newblock {\em The Physics of Fluids}, 24(9):1730--1749, 1981.

\bibitem{CaBr2009}
J.~R. Cary and A.~J. Brizard.
\newblock Hamiltonian theory of guiding-center motion.
\newblock {\em Rev. Mod. Phys.}, 81:693--738, 2009.

\bibitem{Arnold_1989}
V.~I. Arnold.
\newblock {\em Mathematical Methods of Classical Mechanics}.
\newblock Springer, 1989.

\bibitem{Cary_1983}
J.~R. Cary and R.~G. Littlejohn.
\newblock Noncanonical {H}amiltonian mechanics and its application to magnetic
  field line flow.
\newblock {\em Ann. Phys.}, 151:1, 1983.

\bibitem{Qin_2008}
H.~Qin and X.~Guan.
\newblock Variational symplectic integrator for long-time simulations of the
  guiding-center motion of charged particles in general magnetic fields.
\newblock {\em Phys. Rev. Lett.}, 100(3):035006, 2008.

\bibitem{Ellison_thesis}
C.~Leland Ellison.
\newblock {\em Development of Multistep and Degenerate Variational Integrators
  for Applications in Plasma Physics}.
\newblock Doctoral {T}hesis, Princeton University, 2016.

\bibitem{Ellison_2018}
C.~L. Ellison, J.~M. Finn, J.~W. Burby, M.~Kraus, H.~Qin, and W.~M. Tang.
\newblock Degenerate variatonal integrators for magnetic field line flow and
  guiding center trajectories.
\newblock {\em Physics of Plasmas}, 25:052502, 2017.

\bibitem{Hairer_2006}
E.~Hairer, C.~Lubich, and G.~Wanner.
\newblock {\em Geometric Numerical Integration}.
\newblock Springer, 2006.

\bibitem{Rowley_2002}
C.~W. Rowley and J.~E. Marsden.
\newblock Variational integrators for degenerate {L}agrangians, with
  application to point vortices.
\newblock {\em 41st IEEE Conference on Decision and Control}, 40:1521, 2002.

\bibitem{Tao_2016}
M.~Tao.
\newblock {Explicit symplectic approximation of nonseparable Hamiltonians:
  Algorithm and long-time performance}.
\newblock {\em Phys. Rev. E}, 94:043303, 2016.

\end{thebibliography}

\end{document}